\documentclass[10pt,a4paper,draft]{article}

\title{Convergence analysis of nonconform $H(\div)$-finite elements for the damped time-harmonic Galbrun's equation%
\thanks{
The author acknowledges support from DFG project 468728622 and DFG SFB 1456 project 432680300. 
}}
\author{Martin Halla\thanks{
Institut f\"ur Numerische und Angewandte Mathematik,
Georg-Augst Universität Göttingen,
Lotzestr.\ 16-18, 37083 Göttingen, Deutschland.
e-mail: m.halla@math.uni-goettingen.de
}
}

\usepackage{a4wide}
\usepackage{amsmath,amsthm,amsfonts,amssymb}
\usepackage{mathtools}
\usepackage{enumitem}\setitemize{leftmargin=5mm}
\usepackage{accents}
\usepackage{color}
\usepackage[usenames,dvipsnames,svgnames,table]{xcolor}
\usepackage{cite}
\usepackage{soul}
\usepackage[noabbrev,nameinlink]{cleveref}
\usepackage{fancyhdr}
\usepackage{xfrac}
\usepackage[makeroom]{cancel}
\theoremstyle{plain}
\newtheorem{theorem}{Theorem}
\newtheorem{lemma}[theorem]{Lemma}

\theoremstyle{definition}

\newtheorem{remark}[theorem]{Remark}

\newcommand{\ol}[1]{\overline{#1}}
\newcommand{\ull}[1]{\underline{#1}}
\newcommand{\spl}{\langle}
\newcommand{\spr}{\rangle}
\newcommand{\bpm}{\begin{pmatrix}}
\newcommand{\epm}{\end{pmatrix}}

\newcommand{\bb}{(}

\newcommand{\bjump}[1]{[\![ #1 ]\!]_{\bflow}}

\newcommand{\avg}[1]{\{\!\!\{#1\}\!\!\}}
\newcommand{\h}{\mathfrak{h}}
\newcommand{\hf}{h_\face}
\newcommand{\htet}{h_\tet}

\newcommand{\hmax}{h_n} %{\overline{\mathfrak{h}}_n}
\newcommand{\ALpen}{\alpha}

\newcommand{\dn}{d_n} %{\mathrm{d}_n}
\newcommand{\Pk}{P_k}

\newcommand{\ses}{a}
\newcommand{\sesn}{a_n}

\newcommand{\Find}{{\calF_n^\mathrm{int}}}

\newcommand{\dom}{\mathcal{O}} % main domain
\newcommand{\nv}{\boldsymbol{\nu}}
\newcommand{\bflow}{\mathbf{b}}
\newcommand{\vel}{\mathbf{b}}
\newcommand{\csu}{\ol{c_s}}
\newcommand{\csl}{\ull{{c_s}}}
\newcommand{\rhou}{\ol{\rho}}
\newcommand{\rhol}{\ull{\rho}}
\newcommand{\gammal}{\ull{\gamma}}
\newcommand{\gammau}{\ol{\gamma}}
\newcommand{\angvel}{\Omega}
\newcommand{\opd}{(\omega+i\conv+i\angvel\times)}
\newcommand{\opdn}{(\omega+i\Db+i\angvel\times)}
\renewcommand{\u}{\mathbf{u}}
\renewcommand{\v}{\mathbf{v}}
\newcommand{\w}{\mathbf{w}}
\newcommand{\z}{\mathbf{z}}
\newcommand{\m}{\underline{\underline{m}}} % short notation for Section 3.4
\newcommand{\q}{\mathbf{q}} % scaled pressure gradient

\newcommand{\nh}{n} % index for Section 2
\newcommand{\limnh}{n\to\infty} %¸ limit for Section 2
\newcommand{\Bstab}{B}
\newcommand{\Bcomp}{K}

\newcommand{\vsi}{v}
\newcommand{\vsii}{\tilde v}
\newcommand{\vi}{\v}
\newcommand{\vii}{\tilde \v}
\newcommand{\wi}{\w}

\newcommand{\Htbc}{H^2_{*,\mathrm{Neu}}}
\newcommand{\Hz}{\bH^1_{\nv0}}
\newcommand{\Hpw}{\bH^1_\mathrm{pw}(\calT_n)}

\newcommand{\X}{\IX}
\newcommand{\Xn}{\X_n}

\newcommand{\Qn}{Q_n}
\newcommand{\bQn}{\mathbf{Q}_n}

\newcommand{\tet}{\tau}
\newcommand{\tetp}{D_\tet} %patch
\newcommand{\tetf}{D_\face} %patch
\newcommand{\face}{F}

\newcommand{\Cpis}{C_\pi^\#}
\newcommand{\Cpih}{\tilde C_\pi}
\newcommand{\CKeta}{C_1}
\newcommand{\CtKn}{C_2}

\newcommand{\Ctr}{C_\mathrm{dt}}
\newcommand{\Chratio}{C_\mathrm{sh}}
\newcommand{\Capr}{C_\mathrm{apr}}
\newcommand{\Cface}{C_\mathrm{ab}}
\newcommand{\Cinv}{C_\mathrm{inv}}
\newcommand{\CL}{C^L}
\newcommand{\Cm}{C_{\m}}

\newcommand{\CtB}{C_{\tilde B}}
\newcommand{\CYi}{C_{Y,1}}
\newcommand{\CYii}{C_{Y,2}}
\newcommand{\Mach}{\|c_s^{-1}\bflow\|_{\bL^\infty}}

\newcommand{\g}{\mathbf{g}}
\newcommand{\bpsi}{\boldsymbol{\psi}}

\newcommand{\tf}{\bpsi_n}
\newcommand{\LO}{R_n}
\newcommand{\Db}{D^n_\vel}
\newcommand{\Kuni}{M} % finite dimensional operator used for the def. of v_n
\newcommand{\Keta}{K_G} % compact operator from HallaHohage:21
\newcommand{\PV}{P_V}
\newcommand{\tPV}{P_{\tilde V_n}}
\newcommand{\PVn}{P_{V_n}}
\newcommand{\PLz}{P_{L^2_0}}
\newcommand{\On}{O_n} % operator which tends to zero
\newcommand{\tOn}{\tilde O_n} % operator which tends to zero
\newcommand{\hOn}{\hat O_n} % operator which tends to zero
\newcommand{\cOn}{\check O_n} % operator which tends to zero
\renewcommand{\Re}{\operatorname{Re}}
\renewcommand{\Im}{\operatorname{Im}}
\newcommand{\conv}{\partial_{\bflow}} %{(\mathbf{\bflow}\cdot \nabla)}
\renewcommand{\div}{\operatorname{div}}

\DeclareMathOperator{\hess}{Hess}

\DeclareMathOperator{\Id}{Id}

\DeclareMathOperator{\supp}{supp}

\DeclareMathOperator{\dist}{dist}

\DeclareMathOperator{\sign}{sgn}

\DeclareMathOperator{\mean}{mean}

\newcommand{\inner}[1]{\langle #1 \rangle}

\newcommand\restr[2]{{ \left.\kern-\nulldelimiterspace #1 \vphantom{\big|} \right|_{#2} }}

\newcommand{\IN}{\mathbb{N}}

\newcommand{\IR}{\mathbb{R}}

\newcommand{\IX}{\mathbb{X}}

\newcommand{\calF}{\mathcal{F}}

\newcommand{\calJ}{\mathcal{J}}
\newcommand{\calK}{\mathcal{K}}

\newcommand{\calT}{\mathcal{T}}
\newcommand{\calU}{\mathcal{U}}

\renewcommand{\bf}{\mathbf{f}}
\newcommand{\bg}{\mathbf{g}}

\newcommand{\bu}{\mathbf{u}}

\newcommand{\bx}{\mathbf{x}}

\newcommand{\bC}{\mathbf{C}}

\newcommand{\bH}{\mathbf{H}}

\newcommand{\bL}{\mathbf{L}}

\newcommand{\bP}{\mathbf{P}}
\newcommand{\bQ}{\mathbf{Q}}

\newcommand{\bV}{\mathbf{V}}
\newcommand{\bW}{\mathbf{W}}
\newcommand{\bX}{\mathbf{X}}

\definecolor{pscol}{rgb}{0.8,0,0}

\usepackage{comment} 
\excludecomment{extproof}

\begin{document}

\maketitle
\begin{abstract}
\noindent
We consider the damped time-harmonic Galbrun's equation, which is used to model stellar oscillations.
We introduce a discontinuous Galerkin finite element method (DGFEM) with $H(\div)$-elements, which is nonconform with respect to the convection operator.
We report a convergence analysis, which is based on the frameworks of discrete approximation schemes and T-compatibility.
A novelty is that we show how to interprete a DGFEM as a discrete approximation scheme and this approach enables us to apply compact perturbation arguments in a DG-setting, and to circumvent any extra regularity assumptions on the solution.
The advantage of the proposed $H(\div)$-DGFEM compared to $H^1$-conforming methods is that we do not require a minimal polynomial order or any special assumptions on the mesh structure.
The considered DGFEM is constructed without a stabilization term, which considerably improves the assumption on the smallness of the Mach number compared to other DG methods and $H^1$-conforming methods, and the obtained bound is fairly explicit.
In addition, the method is robust with respect to the drastic changes of magnitude of the density and sound speed, which occur in stars.
The convergence of the method is obtained without additional regularity assumptions on the solution, and for smooth solutions and parameters convergence rates are derived.
% {\color{blue}Computational experiments to accompany the theoretical results are included}.
\vspace{2mm}\\
\textbf{Key words. }Galbrun's equation, stellar oscillations, T-compatibility, T-coercivity,  discrete approximation scheme, discontinuous Galerkin, DGFEM
\vspace{1mm}\\
\textbf{MSC codes. }35L05, 35Q85, 65N30
% 35L05 Wave equation
% 35Q35 PDEs in connection with fluid mechanics
% 35Q85 PDEs in connection with astronomy and astrophysics
% 85A20 Planetary atmospheres 
% 65N30 Finite element, Rayleigh-Ritz and Galerkin methods for boundary value problems involving PDEs 
\end{abstract}
\sloppy
% \tableofcontents

\section{Introduction}

In this article we introduce, analyze and test a particular finite element method to approximatively solve the damped time-harmonic Galbrun's equation
\begin{subequations}\label{eq:Galbrun}
\begin{align}
\begin{aligned}
- \nabla\left(\rho c_s^2\div \u\right)
+(\div \u) \nabla p
&-\nabla(\nabla p\cdot \u) 
-\rho\opd^2\u \\
&+(\hess(p)-\rho\hess(\phi))\u 
+ \gamma \rho (-i \omega) \u
= \bf \quad \mbox{in } \dom,
\end{aligned}\\
&\hspace{-25mm}\nv\cdot\bu=0\quad \mbox{on } \partial\dom,
\end{align}
\end{subequations}
where $\rho, p, \phi, c_s,\bflow, \angvel$ and $\bf$ denote density, pressure, gravitational potential, sound speed, background velocity, angular velocity of the frame and sources, $\partial_\bflow := \sum_{l=1}^3 \bflow_l\partial_{\bx_l}$ denotes the directional derivative in direction $\bflow$, $\hess(p)$ the Hessian of $p$, $\dom\subset\mathbb{R}^3$ a bounded domain, and damping is modeled by the term $- i \omega \gamma \rho \u$ with damping coefficient $\gamma$.
The Galbrun's equation was first derived in \cite{Galbrun:31} as a linearization of the nonlinear Euler's equation and serves as a model in aeroacoustics \cite{maeder202090} and in an extended form in asterophysics \cite{LyndenBOstriker:67}.
In the time-domain the Galbrun's equation was analyzed in \cite{HaeggBerggren:19}.
In the time-harmonic domain a well-posedness analysis in an aeroacoustic setting was reported in \cite{BonnetBDMercierMillotPernetPeynaud:12} through the introduction of an additional transport equation.
Different to that in a stellar context well-posedness results were reported in \cite{HallaHohage:21,Halla:22GalExt} by exploiting the damping effects in stars.
Concerning the numerical approximation of Galbrun's equation it is well known \cite{chabassier:hal-01833043} that naive discretizations may yield unreliable results.
To construct stable methods a path is to follow the $T$-analysis from the continuous level \cite{HallaHohage:21} and to try to mimic the analysis on the discrete level.
However, the $T$-compatibility framework \cite{Halla:21Tcomp}, which was previously applied successfully to various Maxwell problems \cite{Halla:21SteklovAppr,Unger:21,Halla:21SC} and perfectly matched layer methods \cite{Halla:21PML,Halla:22PMLani}, turned out to require too strong assumptions to allow the analysis of discretizations to the Galbrun's equation.
As a remedy a version with weaker assumptions was introduced in \cite{HallaLehrenfeldStocker:22} and successfully applied for the convergence analysis of \emph{divergence stable} $\bH^1$-conforming finite element discretizations of Galbrun's equation.
Therein the so-called \emph{divergence stability} is ensured by assumptions on the mesh structure and the polynomial order of the method.
However, those methods require a lot of computational cost, e.g.\ in three dimensions general meshes are speculated to require a minimal polynomial degree between six and eight \cite{Zhang:09,Zhang:11a} and barycentric refined meshes require a minimal polynomial degree three \cite{GuzmanNeilan:18}.
Although note that the barycentric refinement produces a lot additional degrees of freedom without a reduction of the element diameters.

In this article we employ $H(\div)$-conforming finite elements and treat the nonconformity with respect to the convection operator with a discontinuous Galerkin technique.
In particular, we apply a reconstruction operator to lift the jumps and avoid a stablization term to optimize the assumption on the smallness of the Mach number.
The obtained method does not require any special assumption on the meshes and works for all polynomials orders greater equal than one.
Most importantly the method is robust with respect to the drastic changes of magnitude of the density and sound speed, which occur in stars.
Apart from proving convergence of the method a major contribution of this article is to show how to interprete a DGFEM as a discrete approximation scheme.

The remainder of this article is structured as follows.
In \cref{sec:dddfem} we introduce the applied $H(\div)$-discontinuous Galerkin finite element method.
In \cref{sec:framework} we recall the abstract framework from \cite{HallaLehrenfeldStocker:22} and show that our $H(\div)$-DGFE-method constitutes an asymptoticly consistent \emph{discrete approximation scheme}.
In \cref{sec:convergence} we introduce discrete operators $(T_n)_{n\in\IN}$, analyze their properties and report our main convergence result in \cref{thm:convergence}.
% In \cref{sec:numex} we present computational examples.

\section{Formulation of the $H(\div)$-DDG-FEM}\label{sec:dddfem}

For two Hilbert spaces $(X, \inner{\cdot,\cdot}_X)$, $(Y, \inner{\cdot,\cdot}_Y)$ let $L(X,Y)$ be the space of bounded linear operators from $X$ to $Y$, and set $L(X):=L(X,X)$.
For any space $X$ of scalar valued functions let $\bX:=X^3$.
For $A\in L(X)$ let the bounded sesquilinear form $a(\cdot,\cdot)$ be defined by the relation
\begin{align}\label{eq:RelAa}
\spl Au,u'\spr_X=a(u,u') \quad\text{for all }u,u'\in X,
\end{align}
and vice-versa for a given bounded sesquilinear form $a(\cdot,\cdot)$ let $A\in L(X)$ be defined by the relation \eqref{eq:RelAa}.
We call an operator $A\in L(X)$ coercive, if $\sup_{u\in X\setminus\{0\}} |\spl Au,u\spr|/\|u\|_X^2>0$.
For a bijective operator $T\in L(X)$ we call $A\in L(X)$ to be weakly left (right) $T$-coercive, if there exists a compact operator $K\in L(X)$ such that $T^*A+K$ ($AT+K$) is coercive.
Note that historically the notion of left $T$-coercivity was used, because thence the operator $T$ selects a suitable test function: $\spl T^*Au,u\spr_X=a(u,Tu)$.
However, when conducting the stability/compatibility analysis on the discrete level the notion of right weak $T$-coercivity seems favorable \cite{HallaLehrenfeldStocker:22}, because it avoids the introduction and subsequent treatment of adjoint operators.
For expressions $A,B\in\IR$ we employ the notation $A\lesssim B$, if there exists a constant $C>0$ such that $A\leq CB$.
The constant $C>0$ may be different at each occurence and can depend on the domain $\dom$, the physical parameters $\rho,c_s,p,,\phi,\gamma,\bflow,\omega,\angvel$, and on the sequence of Galerkin spaces $(X_n)_{n\in\IN}$.
However, it will always be independent of the index $n$ and any involved functions ($u,v\in X, u_n\in X_n$, etc.) which may appear in the terms $A$ and $B$.

\subsection{Variational formulation}

Let $\dom\subset\IR^3$ be a bounded convex Lipschitz polyhedron, $\omega\in\IR\setminus\{0\}$ and $\angvel\in\IR^3$.
For brevity all function spaces without specified domain are considered on the domain $\dom$, e.g.\ $L^2=L^2(\dom)$, etc..
Let $c_s,\rho\in W^{1,\infty}(\dom,\IR)$, $\gamma\in L^\infty(\dom,\IR)$ and constants $\csl, \csu, \rhol, \rhou, \gammal,\gammau>0$ be such that $\csl\leq c_s(\bx)\leq \csu$, $\rhol \leq \rho(\bx)\leq \rhou$ and $\gammal\leq \gamma(\bx)\leq \gammau$ for all $\bx\in\dom$.
In addition let $\vel\in W^{1,\infty}(\dom,\IR^3)$ be compactly supported in $\dom$ and $p,\phi \in W^{2,\infty}(\dom,\IR)$´.
For a scalar function $u$ we consider its gradient to be a column vector $\nabla u:=(\partial_{\bx_1}u,\partial_{\bx_2}u,\partial_{\bx_3}u)^\top$, its Hessian to be a matrix $\hess(u):=(\partial_{\bx_n}\partial_{\bx_m}u)_{n,m=1,2,3}$, and for a (column) vectorial function $\u$ we consider its gradient to be a matrix $\nabla\u:=(\partial_{\bx_m}\u_n)_{n,m=1,2,3}$.
We abreviate the $L^2$- and $\bL^2$-scalar products as $\spl\cdot,\cdot\spr$.
For any space $X\subset L^2$ let $X_*:=\{u\in X\colon \spl u,1\spr=1\}$ and $L^2_0:=L^2_*$.
For functions $f\in W^{1,\infty}, \bf\in \bW^{1,\infty}$ we denote their Lipschitz constants as $\CL_f$ and $\CL_{\bf}$.
In addition we introduce the space $\Hz:=\{\u\in\bH^1\colon \nv\cdot\u=0\text{ on }\partial\dom\}$ and the weighted semi norm $|\u|_{\bH^1_{c_s^2\rho}}^2:=\|c_s\rho^{\sfrac12}\nabla\u\|_{(L^2)^{3\times3}}^2$.
Let $\conv:=\vel\cdot\nabla=\sum_{l=1}^3 \vel_l\partial_{\bx_l}$.
Thence we introduce
\begin{align*}
\IX&:=\{\u\in \bL^2\colon \div\u\in L^2,\quad \conv \u\in \bL^2,\quad \nv\cdot\u=0 \text{ on } \partial \dom\},\\
\spl\u,\u'\spr_\IX&:=\inner{\div \u,\div \u'} + \inner{\conv \u,\conv \u'} + \inner{\u,\u'},
\end{align*}
which constitutes a Hilbert space \cite[Lemma~2.1]{HallaHohage:21}.
We recall that $\bC_0^\infty$ is dense in $\IX$ \cite[Theorem~6]{HallaLehrenfeldStocker:22}.
Let
\begin{align}\label{eq:ses}
\begin{aligned}
\ses(\u,\u'):=&\,\inner{c_s^2\rho\div \u,\div \u'} - \inner{\rho\opd \u,\opd \u'}\\
&+\inner{\div\u,\nabla p\cdot \u'}+\inner{\nabla p\cdot\u,\div \u'}+\inner{(\hess(p)-\rho\hess(\phi)) \u,\u'} \\
&- i\omega\inner{\gamma\rho\u,\u'}
\end{aligned}
\end{align}
for all $\u,\u'\in\IX$.
Then assuming the conservation of mass $\div(\rho\vel)=0$ the variational formulation of \eqref{eq:Galbrun} is $\ses(\u,\u')=\spl\bf,\u'\spr$ for all $\u'\in\IX$ \cite[Section~2.3]{HallaHohage:21}.

\subsection{$H(\div)$-finite elements}

Let $(\calT_n)_{n\in\IN}$ be a sequence of tetrahedral meshes of $\dom$, and $\calF_n$ $ (\calF_n^\mathrm{int})$ be the collection of the (interior) faces of $\calT_n$.
For $\tet\in\calT_n$ denote $\calF_\tet\subset\calF_n$ the faces of $\tet$, and $\tetp:=\mathrm{int}\,\ol{\{\bigcup\tilde\tet\colon \tilde\tet\in\calT_n, \ol{\tilde\tet}\cap\ol{\tet}\neq\emptyset\}}$ the macro element of $\tet$ which consists of all elements which share a face with $\tet$.
Also for $\face\in\calF_n$ we define the macro element $\tetf:=\mathrm{int}\,\ol{\{\bigcup\tet\colon \face\in\calF_{\tet}, \tet\in\calT_n\}}$.
Further for $\calU\in\{\tet,\tetp,\face,\tetf\colon \tet\in\calT_n,\face\in\calF_n\}$ and a scalar function $g$ we abreviate $\ull{g}_\calU:=\inf_{\bx\in\calU}g(\bx)$ and $\ol{g}_\calU:=\sup_{\bx\in\calU}g(\bx)$.
For $\tet\in\calT_n$ and $\face\in\calF_n$ let $\htet$ and $\hf$ respectively be their diameters.
Let $\h\colon\calF_n\to\IR$ be defined by $\h|_\face:=\hf$ for $\face\in\calF_n$ and $\hmax:=\max_{\tet\in\calT_n}\htet$.
We assume that $\lim_{n\to\infty}\hmax=0$ and that $(\calT_n)_{n\in\IN}$ is shape regular, i.e. there exists a constant $\Chratio>0$ such that
\begin{align}\label{eq:Chratio}
\htet\leq\Chratio\hf \quad\text{for all }\tet\in\calT_n, \face\in\calF_\tet, n\in\IN.
\end{align}
Let $\Pk$ be the space of scalar polynomials with maximal degree $k\in\IN_0$.
In this article we consider Brezzi–Douglas–Marini elements \cite[Chapter~14.5.1]{ErnGuermond:FEi} due to their approximations properties.
In principle other $H(\div)$-elements such as Raviart-Thomas elements can easily be treated with the framework of this article, whereat some details deserve to be treated with care.
% RT: minimal $l\geq1$ necessary? convergence order ...
% RT
% and $\RTk:=(\Pk)^3\oplus \bx\Pk$ be the Raviart-Thomas space ($\bx\in\IR^3$).
% Note that in this article we choose Raviart-Thomas elements only for convenience and the presented theory can be applied to other $H(\div)$-compatible elements e.g.\ such as the Brezzi–Douglas–Marini elements without any changes.
We consider the polynomial degree $k\in\IN=\{1,2,\dots\}$ to be uniform and fixed in the entire article.
We introduce the finite element spaces
\begin{align*}
% RT
% \Xn:=\{\u\in H_0(\div)\colon \u|_\tet\in\RTk\text{ for all }\tet\in\calT_n\},\qquad
% \Qn:=\{u\in L^2_0\colon u|_\tet\in\Pk\text{ for all }\tet\in\calT_n\}.
\Xn&:=\{\u\in H_0(\div)\colon \u|_\tet\in\bP_{k}\text{ for all }\tet\in\calT_n\},\qquad
\Qn:=\{u\in L^2_0\colon u|_\tet\in P_{k-1}\text{ for all }\tet\in\calT_n\},\\
\Xn^\mathrm{wbc}&:=\{\u\in H(\div)\colon \u|_\tet\in\bP_{k}\text{ for all }\tet\in\calT_n\},\qquad\hspace{-1.85mm}
\Qn^\mathrm{wbc}:=\{u\in L^2\colon u|_\tet\in P_{k-1}\text{ for all }\tet\in\calT_n\}.
\end{align*}
We consider the space $\Qn$ to be equipped with the standard $L^2$-scalar product, whereas the scalar product on $\Xn$ will be specified lateron.
For each $\tet\in\calT_n$ let
% RT
% $\pi_\tet^d\colon\bH^s(\tet)\to\RTk$, $s>1/2$
$\pi_\tet^d\colon\bH^s(\tet)\to\bP_k$, $s>1/2$
and $\pi_\tet^l\colon L^2(\tet)\to\Pk$ be the respective standard local interpolation operators, and
$\pi_n^d\colon\bH^s\to\Xn^\mathrm{wbc}$, $s>1/2$, $\pi_n^d|_\tet:=\pi_\tet^d, \tet\in\calT_n$
and
$\pi_n^l\colon L^2\to\Qn^\mathrm{wbc}$, $\pi_n^l|_\tet:=\pi_\tet^l, \tet\in\calT_n$ be their global versions.
Note that $\pi_n^d\v\in\Xn$ if $\nv\cdot\v=0$ on $\partial\dom$ and $\pi_n^lv\in\Qn$ if $v\in L^2_0$.
We recall the commutation $\div\pi_n^d=\pi_n^l\div$ and the approximation (and boundedness) properties \cite[Theorem 11.13]{ErnGuermond:FEi}
% \cite[Proposition~2.5.1]{BoffiBrezziFortin:13}
% \cite[Theorem~16.4]{ErnGuermond:FEi}
\begin{align}\label{eq:EstApr}
|\v-\pi_n^d\v|_{\bH^m(\tet)}
\leq \Capr \htet^{r-m} |\v|_{\bH^r(\tet)}, \qquad
|v-\pi_n^lv|_{H^m(\tet)} \leq \Capr \htet^{r} |v|_{H^r(\tet)}
\end{align}
and an approximation result on faces \cite[Remark~12.17]{ErnGuermond:FEi}
\begin{align}\label{eq:EstFace}
\|\v-\pi_n^d\v\|_{\bL^2(\partial\tet)} \leq \Cface \htet^{r-1/2} |\v|_{\bH^r(\tet)}
\end{align}
with constants $\Capr, \Cface>0$
for all $r\in[1,k+1]$, $m\in[0,r]$, $\v\in\bH^r(\tet)$, $v\in H^r(\tet)$, $\tet\in\calT_n$, $n\in\IN$.
In addition, we recall the discrete trace inequality \cite[Lemma~12.8]{ErnGuermond:FEi}
\begin{align}\label{eq:DisTrIneq}
\|u\|_{L^2(\partial\tet)} \leq \Ctr \htet^{-1/2} \|u\|_{L^2(\tet)}
\end{align}
and the discrete inverse inequality \cite[Lemma~12.15]{ErnGuermond:FEi}
\begin{align}\label{eq:DisInv}
|u|_{H^1(\tet)} \leq \Cinv \htet^{-1} \|u\|_{L^2(\tet)}
\end{align}
with constants $\Ctr,\Cinv>0$
% RT
% for all $u\in P_{k+1}$, $\tet\in\calT_n$, $n\in\IN$.
for all $u\in P_{k}$, $\tet\in\calT_n$, $n\in\IN$.

\subsection{Distributional discontinuous Galerkin method}

Let $\Hpw:=\{\u\in\bL^2\colon \u|_\tet\in\bH^1(\tet) \text{ for all }\tet\in\calT_n\}$.
For $\face\in\calF_n^\mathrm{int}$ we denote the neighboring elements of $\face$ as $\tet_1,  \tet_2\in\calT_n$ and set $\tet_\face:=\ol{\tet_1\cup\tet_2}$.
Thence for $\u\in\Hpw$ and $F\in\calF_n^\mathrm{int}$ we denote the traces of $\bu|_{\tet_1}$ and $\bu|_{\tet_2}$ on $\face$ as $\u_1$ and $\u_2$ respectively.
Thus for $\u\in\Hpw$ and $\face\in\calF_n^\mathrm{int}$ we define the following average and jump terms
\begin{align*}
\avg{\u}:=\frac12(\u_1+\u_2),
&& \bjump{\u}:=(\vel\cdot\nv_1)\u_1+(\vel\cdot\nv_2)\u_2.
\end{align*}
In addition we introduce the abbreviations
\begin{align*}
\spl \cdot,\cdot \spr_\Find:=
\sum_{\face\in\calF_n^\mathrm{int}} \spl \cdot,\cdot \spr_{\bL^2(\face)},\quad
\| \cdot \|_\Find^2:=\spl \cdot,\cdot \spr_\Find.
\end{align*}
We introduce a lifting operator (see, e.g., \cite[Chapter~4.3]{DiPietroErn:12}) related to the differential operator $\conv$.
Note that we choose a sign convention for $\LO^\face$ as in \cite{BuffaOrtner:09}, which is opposite to the one used in \cite[Chapter~4.3]{DiPietroErn:12}.
% RT
% Let $l\in \{0,\dots,k+1\}$ and
Let $l\in\IN$.
\begin{align*}
\bQn:=\{\tf\in\bL^2\colon \tf|_\tet\in \bP_{l} \text{ for all }\tet\in\calT_n\}.
\end{align*}
Then for $\u_n\in\Xn$ and $\face\in\calF_n^\mathrm{int}$ let $\LO^\face\u_n\in\bQn$ be the solution to
\begin{align*}
\spl \LO^\face\u_n,\tf \spr=-\spl \bjump{\u_n}, \avg{\tf} \spr_{\bL^2(\face)} \quad\text{for all }\tf\in\bQn.
\end{align*}
We observe that $\supp(\LO^\face\u_n)=\tet_\face$ and it easily follows with \eqref{eq:DisTrIneq} that
\begin{align}\label{eq:LObound}
\|\LO^\face\u_n\|_{\bL^2}\leq \Ctr \|\h^{-1/2}\bjump{\u_n}\|_{\bL^2(\face)}
\quad\text{for all }\u_n\in\Xn, \face\in\calF_n^\mathrm{int}, n\in\IN.
\end{align}
Thence we define $\LO:=\sum_{\face\in\calF_n^\mathrm{int}}\LO^\face$ and the linear operator $\Db\colon \Xn\to\bQn$ by
\begin{align*}
(\Db\u_n)|_{\tet}:=\conv(\u_n|_{\tet})+\LO\u_n \quad\text{for all }\tet\in\calT_n.
\end{align*}
We remark that for less complicated equations the natural choice for $l$ is $l=k-1$, while $l=k$ might simplify the implementation \cite{DiPietroErn:12}.
However, for the DG method applied in this article it is advisable to choose indeed $l=k$ to exploit the full potential convergence rate (see \cref{thm:convergence}).
Also note that $l\geq1$ (i.e.\ $l=0$ is excluded) is necessary to obtain \cref{lem:weakcompactness} and subsequent results.
We introduce the following scalar product on $\Xn$
\begin{align*}
\spl \u_n,\u_n' \spr_{\Xn}:=\spl \div\u_n,\div\u_n' \spr
+\spl \u_n,\u_n' \spr
+\spl \Db\u_n,\Db\u_n' \spr.
\end{align*}
Now let
\begin{align}\label{eq:sesn}
\begin{aligned}
\sesn(\u_n,\u_n'):=&\,\inner{c_s^2\rho\div \u_n,\div \u'_n} - \inner{\rho\opdn \u_n,\opdn \u_n'}\\
&+\inner{\div\u_n,\nabla p\cdot \u_n'}+\inner{\nabla p\cdot\u_n,\div \u_n'}+\inner{(\hess(p)-\rho\hess(\phi)) \u_n,\u_n'} \\
&- i\omega\inner{\gamma\rho\u_n,\u_n'}
\quad\text{for all }\u_n,\u_n'\in\Xn.
\end{aligned}
\end{align}
Thus it easily follows that $\sup_{n\in\IN}\|A_n\|_{L(\Xn)}<\infty$.

\section{Abstract framework}\label{sec:framework}

\subsection{Discrete approximation schemes and T-compatibility}

We remark that in this section we use the symbols $\tilde A, \tilde A_n$ for generic operators, because the symbols $A, A_n$ are already occupied due to the introduction of the sesquilinears $\ses(\cdot,\cdot)$ and $\sesn(\cdot,\cdot)$ in \eqref{eq:ses} and \eqref{eq:sesn}.
For a Hilbert space $X$ and $\tilde A\in L(X)$ we consider \emph{discrete approximation schemes} (DAS) of $(X,\tilde A)$ in the following way.
Let $(X_n)_{n\in\IN}$ be a sequence of finite dimensional Hilbert spaces and $\tilde A_n\in L(X_n)$.
Note that we do not demand that the spaces $X_n$ are subspaces of $X$. 
Instead we demand that there exist operators $p_n \in L(X,X_n)$ such that $\lim_{n\to\infty} \|p_n u\|_{X_n}=\|u\|_X$ for each $u\in X$.
We then define the following properties of a discrete approximation scheme:
\begin{itemize}[leftmargin=*]
\setlength\itemsep{0em}
\item A sequence $(u_n)_{n\in\IN}, u_n\in X_n$ is said to \emph{converge} to $u\in X$, if $\lim_{n\to\infty}\|p_nu-u_n\|_{X_n}=0$.
\item A sequence $(u_n)_{n\in\IN}, u_n\in X_n$ is said to be \emph{compact}, if for every subsequence $\IN'\subset\IN$ there exists a subsubsequence $\IN''\subset\IN'$ such that $(u_n)_{n\in\IN''}$ converges (to a $u\in X$).
\item A sequence $(\tilde A_n)_{n\in\IN}, \tilde A_n\in L(X_n)$ is said to be \emph{asymptotic consistent} or to \emph{approximate} $\tilde A\in L(X)$, if $\lim_{n\to\infty}\|\tilde A_np_nu-p_n\tilde Au\|_{X_n}=0$ for each $u\in X$.
\item A sequence of operators $(\tilde A_n)_{n\in\IN}, \tilde A_n\in L(X_n)$ is said to be \emph{compact}, if for every bounded sequence $(u_n)_{n\in\IN}, u_n\in X_n$, $\|u_n\|_{X_n}\leq C$ the sequence $(\tilde A_nu_n)_{n\in\IN}$ is compact.
\item A sequence of operators $(\tilde A_\nh)_{\nh\in\IN}, \tilde A_\nh\in L(X_\nh)$ is said to be \emph{stable}, if there exist constants $C,\nh_0>0$ such that $\tilde A_\nh$ is invertible and $\|\tilde A_\nh^{-1}\|_{L(X_\nh)} \leq C$ for all $\nh>\nh_0$.
\item A sequence of operators $(\tilde A_n)_{n\in\IN}, \tilde A_n\in L(X_n)$ is said to be \emph{regular}, if $\|u_n\|_{X_n}\leq C$ and the compactness of $(\tilde A_nu_n)_{n\in\IN}$ imply the compactness of $(u_n)_{n\in\IN}$.
\end{itemize}
The following theorem will be our tool to prove the regularity of approximations.
\begin{theorem}[Theorem~3 of \cite{HallaLehrenfeldStocker:22}]\label{thm:Tcomp}
Assume a constant $C>0$, sequences $(\tilde A_\nh)_{\nh\in\IN}, (T_\nh)_{\nh\in\IN}, (B_\nh)_{\nh\in\IN}, (K_\nh)_{\nh\in\IN}$ and $B\in L(X)$ which satisfy the following:
for each $\nh\in\IN$ it holds $\tilde A_\nh,T_\nh,B_\nh,K_\nh\in L(X_\nh)$,
$\|T_\nh\|_{L(X_\nh)}, \|T_\nh^{-1}\|_{L(X_\nh)}, \|\Bstab_\nh\|_{L(X_\nh)}, \|\Bstab_\nh^{-1}\|_{L(X_\nh)} \leq C$, $B$ is bijective, $(K_\nh)_{\nh\in\IN}$ is compact,
\begin{align*}
\lim_{\limnh} \|T_\nh p_\nh u-p_\nh Tu\|_{X_\nh}=0
\quad\text{and}\quad
\lim_{\limnh} \|\Bstab_\nh p_\nh u - p_\nh \Bstab u\|_{X_\nh}=0
\text{ for each } u\in X,
\end{align*}
and
\begin{align*}
\tilde A_\nh T_\nh=\Bstab_\nh+\Bcomp_\nh.
\end{align*}
Then $(\tilde A_\nh)_{\nh\in\IN}$ is regular.
\end{theorem}
We recall the following lemma, which shows that regularity implies convergence.
\begin{lemma}[Lemmas 1 \& 2 of \cite{HallaLehrenfeldStocker:22}]\label{lem:DASstable}
Let $\tilde A\in L(X)$ be bijective and $(X_n,\tilde A_n,p_n)_{n\in\IN}$ be a discrete approximation scheme of $(X,\tilde A)$ which is regular and asymptotic consistent.
Then $(\tilde A_n)_{n\in\IN}$ is stable.
If $u,u_n$ are the solutions to $\tilde Au=f$ and $\tilde A_n u_n=f_n\in X_n$, and $\lim_{n\to\infty} \|p_nf-f_n\|_{X_n}=0$, then $\lim_{n\to\infty} \|p_nu-u_n\|_{X_n}=0$.
\end{lemma}

\subsection{Interpretation of the $H(\div)$-DDG-FEM as DAS}\label{subsec:DG-DAS}

For $\u\in\IX$ let $p_n\u\in\Xn$ be the solution to
\begin{align*}
\spl p_n\u,\u_n'\spr_{\Xn}=\spl \div\u,\div\u_n'\spr_{L^2}
+\spl \u,\u_n'\spr_{\bL^2}
+\spl \conv\u,\Db\u_n'\spr_{\bL^2} \quad\text{for all }\u_n'\in\Xn.
\end{align*}
It easily follows that $p_n\in L(\IX,\Xn)$ and $\|p_n\|_{L(\IX,\Xn)}\leq1$.
In addition, there holds the Galerkin orthogonality
\begin{align}\label{eq:pnGO}
\begin{aligned}
0=\spl \div(\u-p_n\u),\div\u_n'\spr_{L^2}
+\spl \u-p_n\u,\u_n'\spr_{\bL^2}
+\spl \conv\u-\Db p_n\u,\Db\u_n'\spr_{\bL^2}
\end{aligned}
\end{align}
for all $\u_n'\in\Xn$.
In order to analyze $p_n$ further we introduce the distance function $\dn(\u,\u_n)$ between $\u\in\IX$ and $\u_n\in\Xn$ as
\begin{align*}
\dn(\u,\u_n)^2&:=\|\div\u-\div \u_n\|_{L^2}^2+\|\u-\u_n\|_{\bL^2}^2
+\|\conv\u-\Db \u_n\|_{\bL^2}^2.
\end{align*}
The introduction of $\dn(\cdot,\cdot)$ is necessary, because in general the jump $\bjump{\u}$ is not well-defined for $\u\in\IX$.
It can easily be seen that $\dn(\cdot,\cdot)$ satisfies the triangle inequalities
\begin{align*}
\dn(\u,\u_n)\leq \dn(\tilde\u,\u_n)+\|\u-\tilde\u\|_{\IX},\qquad
\dn(\u,\u_n)\leq \dn(\u,\tilde\u_n)+\|\u_n-\tilde\u_n\|_{\Xn}
\end{align*}
for all $\u,\tilde\u\in\IX, \u_n\tilde\u_n\in\Xn$.
\begin{lemma}\label{lem:dnpnpi}
For each $\u\in\Hz$ it holds that $\dn(\u,p_n\u)\leq\dn(\u,\pi_n^d\u)$.
\end{lemma}
\begin{proof}
We compute by means of \eqref{eq:pnGO} that
\begin{align*}
\dn(\u,p_n\u)^2&=
\|\div(\u-p_n\u)\|_{L^2}^2+\|\u-p_n\u\|_{\bL^2}^2
+\|\conv\u-\Db p_n\u\|_{\bL^2}^2\\
&=\spl \div(\u- p_n\u),\div( \u-\pi_n^d\u)\spr_{L^2}
+\spl \u-p_n\u,\u-\pi_n^d\u,\spr_{\bL^2}\\
&+\spl \conv\u-\Db p_n\u,\conv\u-\Db \pi_n^d\u\spr_{\bL^2}\\
&\leq \dn(\u,p_n\u) \dn(\u,\pi_n^d\u)
\end{align*}
which proves the claim.
\end{proof}
\begin{lemma}\label{lem:dnpihn}
For each $\u\in\Hz\cap\bH^{1+s}$, $s>0$ it holds that $\dn(\u,\pi_n^d\u)\lesssim \hmax^s\|\u\|_{\bH^{1+s}}$.
\end{lemma}
\begin{proof}
The claim follows from the definition of $\dn(\cdot,\cdot)$,
\begin{align*}
\|\conv\u-\Db\pi_n^d\u\|_{\bL^2(\tet)} \leq \|\conv\u-\conv\pi_n^d\u\|_{\bL^2(\tet)}+\|\LO\pi_n^d\u\|_{\bL^2(\tet)},
\end{align*}
\eqref{eq:EstApr}, \eqref{eq:EstFace}, \eqref{eq:LObound} and $\bjump{\u}=0$.
\end{proof}
\begin{lemma}\label{lem:dnpn}
For each $\u\in\IX$ it holds that $\lim_{n\to\infty}\dn(\u,p_n\u)=0$.
\end{lemma}
\begin{proof}
Since $\bC_0^\infty$ is dense in $\IX$ \cite[Theorem~6]{HallaLehrenfeldStocker:22} for each $\epsilon>0$ we can choose $\tilde\u\in\bC_0^\infty$ such that $\|\u-\tilde\u\|_{\IX}<\epsilon$.
Thence we estimate
\begin{align*}
\dn(\u,p_n\u)
&\leq \dn(\u,p_n\tilde\u)+\|p_n(\u-\tilde\u)\|_{\Xn}\\
&\leq \dn(\tilde\u,p_n\tilde\u)+\|p_n(\u-\tilde\u)\|_{\Xn}+\|\u-\tilde\u\|_{\IX}\\
&\leq \dn(\tilde\u,p_n\tilde\u)+2\epsilon.
\end{align*}
Thus it follows with \cref{lem:dnpnpi,lem:dnpihn} that $\limsup_{n\to\infty}\dn(\u,p_n\u)\leq2\epsilon$.
Since $\epsilon>0$ was arbitrary the claim follows.
\end{proof}
\begin{lemma}\label{lem:dnpi}
For each $\u\in\Hz$ it holds that $\lim_{n\to\infty}\dn(\u,\pi_n^d\u)=0$.
\end{lemma}
\begin{proof}
In principle we proceed as in the proof of \cref{lem:dnpn}.
However, to construct a suitable smooth approximation $\tilde\u$ which respects the boundary condition $\nv\cdot\tilde\u=0$ we need to introduce some technical details.
Let $\calU_1\subset\calU_2\subset\dom$ be such that $\supp\vel\subset\calU_1$, $\dist(\supp\vel,\partial\calU_1)>0$, $\dist(\calU_1,\partial\calU_2)>0$ and $\dist(\calU_2,\partial\dom)>0$.
Note that such $\calU_1,\calU_2$ exist, because by assumption $\dist(\partial\dom,\supp\vel)>0$.
Then let $\chi\in C^\infty$ be such that $\chi=0$ on $\partial\dom$ and $\chi=1$ in $\calU_2$.
For $\u\in\Hz$ let $\u_1:=\chi\u$ and $\u_2:=(1-\chi)\u$.
Let $\epsilon>0$ be given.
Since $\u_1\in\bH^1_0$ and $\bC^\infty_0$ is dense in $\bH^1_0$ we can find $\tilde\u_1\in\bC_0^\infty$ such that $\|\u_1-\tilde\u_1\|_{\bH^1}<\epsilon$.
Thence we estimate
\begin{align*}
\dn(\u,\pi_n^d\u) &\leq \dn(\u_1,\pi_n^d\u_1)+\dn(\u_2,\pi_n^d\u_2)
\end{align*}
and
\begin{align*}
\dn(\u_1,\pi_n^d\u_1)
&\leq \dn(\tilde\u_1,\pi_n^d\tilde\u_1)+\|\u_1-\tilde\u_1\|_{\IX}+\|\pi_n^d(\u_1-\tilde\u_1)\|_{\Xn}\\
&\lesssim \dn(\tilde\u_1,\pi_n^d\tilde\u_1)+(1+\sup_{m\in\IN}\|\pi_m^d\|_{L(\Hz,\IX_m)})\|\u_1-\tilde\u_1\|_{\bH^1}\\
&\lesssim \hmax \|\tilde\u_1\|_{\bH^2}+\epsilon.
\end{align*}
Thus $\limsup_{n\in\IN} \dn(\u_1,\pi_n^d\u_1)\lesssim\epsilon$ and hence $\lim_{n\in\IN} \dn(\u_1,\pi_n^d\u_1)=0$.
For $\u_2$ we use the smoothing operators $\calK_{\epsilon,0}^d, \calK_{\epsilon,0}^b$ defined in \cite[\bb4.1)]{ErnGuermond:16}, set $\tilde\u_2:=\calK_{\epsilon,0}^d\u_2$ and estimate
\begin{align*}
\dn(\u_2,\pi_n^d\u_2) &\leq \dn(\tilde\u_2,\pi_n^d\tilde\u_2)
+\|\u_2-\tilde\u_2\|_{\IX}+\|\pi_n^d(\u_2-\tilde\u_2)\|_{\Xn}.
\end{align*}
In addition we compute that
\begin{align*}
\|\pi_n^d(\u_2-\tilde\u_2)\|_{\Xn}^2
=\|\pi_n^d(\u_2-\tilde\u_2)\|_{\bL^2}^2
+\|\div\pi_n^d(\u_2-\tilde\u_2)\|_{\bL^2}^2
+\sum_{\tet\in\calT_n} \|\Db\pi_n^d(\u_2-\tilde\u_2)\|_{\bL^2(\tet)}^2.
\end{align*}
Let $q\in(2,6)$.
We use the bound $\|\pi_n^d\v\|_{\bL^2}\lesssim \|\v\|_{\bL^q}+\|\div\v\|_{L^2}$, $\v\in H(\div)\cap \bL^q$ \cite[Chapter~17.2]{ErnGuermond:FEi}.
It follows that
\begin{align*}
\|\pi_n^d(\u_2-\tilde\u_2)\|_{\bL^2}^2
+\|\div\pi_n^d(\u_2-\tilde\u_2)\|_{\bL^2}^2
&=\|\pi_n^d(\u_2-\tilde\u_2)\|_{\bL^2}^2
+\|\pi_n^l\div(\u_2-\tilde\u_2)\|_{\bL^2}^2\\
&\lesssim \|\u_2-\tilde\u_2\|_{\bL^q}^2
+\|\div\u_2-\div\tilde\u_2\|_{\bL^2}^2\\
&=\|(1-\calK_{\epsilon,0}^d)\u_2\|_{\bL^q}^2
+\|(1-\calK_{\epsilon,0}^b)\div\u_2\|_{\bL^2}^2.
\end{align*}
In addition, there exists $\epsilon_0>0$ such that $\tilde\u_2|_{\calU_1}=(\calK_{\epsilon,0}\u_2)|_{\calU_1}=0$ for all $\epsilon\in(0,\epsilon_0)$, which we assume henceforth.
There also exists an index $n_0>0$ such that $\tetp\subset\calU_1$ for all $\tet\in\calT_n$ with $\tet\cap\supp\vel\neq\emptyset$ for all $n>n_0$, which we assume henceforth.
Thus
\begin{align*}
\sum_{\tet\in\calT_n} \|\Db\pi_n^d(\u_2-\tilde\u_2)\|_{\bL^2(\tet)}^2
\leq \sum_{\tet\in\calT_n} \|\vel\|_{\bL^\infty(\tetp)}^2 \|\u_2-\tilde\u_2\|_{\bH^1(\tetp)}^2=0.
\end{align*}
Further we note that $\|\u_2-\tilde\u_2\|_{\IX}=\|\u_2-\tilde\u_2\|_{H(\div)}$, because $\supp\vel\cap\sup(\u_2-\tilde\u_2)=\emptyset$.
Altogether we obtain that
\begin{align*}
\dn(\u_2,\pi_n^d\u_2)
\leq \hmax \|\tilde\u_2\|_{\bH^2}
+\|(1-\calK_{\epsilon,0}^d)\u_2\|_{\bL^q}^2
+\|(1-\calK_{\epsilon,0}^b)\div\u_2\|_{\bL^2}^2
\end{align*}
and hence
\begin{align*}
\limsup_{n\in\IN} \dn(\u_2,\pi_n^d\u_2)
\lesssim \|(1-\calK_{\epsilon,0}^d)\u_2\|_{\bL^q}^2
+\|(1-\calK_{\epsilon,0}^b)\div\u_2\|_{\bL^2}^2.
\end{align*}
Due to the continuous Sobolev embedding $\bH^1\hookrightarrow\bL^q$ it holds that $\u_2\in\bL^q$, and the right hand-side of the former inequality tends to zero for $\epsilon\to0$.
Thus $\lim_{n\in\IN} \dn(\u_2,\pi_n^d\u_2)=0$ and the proof is finished.
\end{proof}
\begin{lemma}\label{lem:pnnorm}
For each $\u\in\IX$ it holds that $\lim_{n\to\infty}\|p_n\u\|_{\Xn}=\|\u\|_{\IX}$.
\end{lemma}
\begin{proof}
We compute
\begin{align*}
\|p_n\u\|_{\Xn}^2
&=\spl p_n\u, p_n\u\spr_{\Xn}
=\spl \div\u,\div p_n\u\spr_{L^2}
+\spl \u,p_n\u\spr_{\bL^2}
+\spl \conv\u,\Db p_n\u\spr_{\bL^2}\\
&=\|\u\|_{\IX}^2
+\spl \div\u,\div(p_n\u-\u)\spr_{L^2}
+\spl \u,p_n\u-\u\spr_{\bL^2}
+\spl \conv\u,\Db p_n\u-\conv\u\spr_{\bL^2}.
\end{align*}
Since
\begin{align*}
|\spl \div\u,\div(p_n\u-\u)\spr_{L^2}
+\spl \u,p_n\u-\u\spr_{\bL^2}
+\spl \conv\u,\Db p_n\u-\conv\u\spr_{\bL^2}|
\leq \|\u\|_{\IX} \dn(\u,p_n \u)
\end{align*}
the claim follows from \cref{lem:dnpn}.
\end{proof}
Thus $(\Xn,A_n,p_n)_{n\in\IN}$ forms a discrete approximation scheme of $(\IX,A)$.
Before we establish the asymptotic consistency of this scheme we need to state a weak compactness result.
\begin{lemma}\label{lem:weakcompactness}
Let $(\u_n)_{n\in\IN}$, $\u_n\in\Xn$ satisfy $\sup_{n\in\IN} \|\u_n\|_{\Xn}<\infty$.
Then there exist $\u\in\IX$ and a subsequence $\IN'\subset\IN$ such that
$\u_n\xrightharpoonup{\bL^2}\u$, $ \div\u_n\xrightharpoonup{L^2}\div\u$ and $ \Db\u_n\xrightharpoonup{\bL^2}\partial_\bflow\u$.
\end{lemma}
\begin{proof}
Since $\u_n,\Db\u_n$ and $\div\u_n$ are bounded sequences in $\bL^2$ and $L^2$ respectively, there exist $\u,\g\in\bL^2$, $q\in L^2$ and a subsequence $\IN'\subset\IN$ such that 
$\u_n\xrightharpoonup{\bL^2}\u$, $\Db\u_n\xrightharpoonup{\bL^2}\g$ and $\div\u_n\xrightharpoonup{L^2}q$.
It remains to show that $\g=\conv\u$ and $q=\div\u$.
To this end we work with a distributional technique.
For conforming differential operators (here the divergence operator) this technique is quite standard and for a reconstructed differential operator it can be found, e.g., in the proof of \cite[Theorem 5.2]{BuffaOrtner:09}.
Let $\bpsi\in\bC^\infty_0$ and $\psi\in C^\infty_0$.
Then
\begin{align*}
\spl q,\psi\spr=\lim_{n\in\IN'} \spl \div\u_n,\psi\spr
=\lim_{n\in\IN'} -\spl \u_n,\nabla\psi\spr
=-\spl \u,\nabla\psi\spr
\end{align*}
and hence $q=\div\u$.
Let $\bpsi_n$ be the lowest order standard $\bH^1$-interpolant of $\bpsi$ on the mesh $\calT_n$.
We compute
\begin{align*}
\spl \Db\u_n,\bpsi\spr&=\spl \Db\u_n,\bpsi-\bpsi_n\spr+\spl \Db\u_n,\bpsi_n\spr\\
&=\spl \Db\u_n,\bpsi-\bpsi_n\spr+\sum_{\tet\in\calT_n} \spl \conv\u_n,\bpsi_n\spr_{\bL^2(\tet)}-\spl \bjump{\u_n},\avg{\bpsi_n}\spr_\Find\\
&=\spl \Db\u_n,\bpsi-\bpsi_n\spr+\sum_{\tet\in\calT_n} \spl \conv\u_n,\bpsi_n\spr_{\bL^2(\tet)}-\spl (\nv\cdot\vel)\u_n,\bpsi_n\spr_{\bL^2(\partial\tet)}\\
&=\spl \Db\u_n,\bpsi-\bpsi_n\spr-\spl \u_n,(\conv+\div(\vel))\bpsi_n\spr\\
&=-\spl \u_n,(\conv+\div(\vel))\bpsi\spr
+\spl \Db\u_n,\bpsi-\bpsi_n\spr+\spl \u_n,(\conv+\div(\vel))(\bpsi-\bpsi_n)\spr.
\end{align*}
Since $\|\bpsi-\bpsi_n\|_{\bH^1}\lesssim \hmax \|\bpsi\|_{\bH^2}$ and $\|\u_n\|_{\Xn}\lesssim1$ it follows that
\begin{align*}
\spl \g,\bpsi\spr
=\lim_{n\to\infty} \spl \Db\u_n,\bpsi\spr
=\lim_{n\to\infty} -\spl \u_n,(\conv+\div(\vel))\bpsi\spr
=-\spl \u,(\conv+\div(\vel))\bpsi\spr
\end{align*}
and hence $\conv\u=\g$.
\end{proof}
\begin{theorem}\label{thm:asympcons}
For each $\u\in\IX$ it holds that $\lim_{n\to\infty} \|A_np_n\u-p_nA\u\|_{\Xn}=0$.
\end{theorem}
\begin{proof}
Let $\u\in\IX$ and $(\u_n)_{n\in\IN}$, $\u_n\in\Xn$, $\|\u_n\|_{\Xn}=1$ be such that
\begin{align*}
\|A_np_n\u-p_nA\u\|_{\Xn} \leq |\spl A_np_n\u-p_nA\u, \u_n \spr_{\Xn}|+1/n.
\end{align*}
For each arbitrary subsequence $\IN'\subset\IN$ we choose $\IN''\subset\IN'$ and $\u'\in\IX$ as in \cref{lem:weakcompactness}.
Then it follows with the definition of $p_n$ that
\begin{align*}
\lim_{n\in\IN''} \spl p_nA\u, \u_n \spr_{\Xn}
&=\lim_{n\in\IN''} \big( \spl \div A\u, \div\u_n \spr+\spl A\u, \u_n \spr+\spl \conv A\u, \Db\u_n \spr\big)\\
&= \spl \div A\u, \div\u' \spr+\spl A\u, \u' \spr+\spl \conv A\u, \conv\u' \spr\\
&=\spl A\u, \u' \spr_{\IX}
=a(\u, \u').
\end{align*}
We further compute that
\begin{align}
\nonumber
\spl &A_n p_n\u, \u_n \spr_{\Xn}
=\sesn(p_n\u, \u_n)\\
\nonumber
&=\inner{c_s^2\rho\div p_n\u,\div \u_n} - \inner{\rho\opdn p_n\u,\opdn \u_n}\\
\nonumber
&+\inner{\div p_n\u,\nabla p\cdot \u_n}+\inner{\nabla p\cdot p_n\u,\div \u_n}+\inner{(\hess(p)-\rho\hess(\phi)) p_n\u,\u_n} \\
\nonumber
&- i\omega\inner{\gamma\rho p_n\u,\u_n}\\
&\hspace{-2mm}\begin{aligned}\label{eq:lima}
\left.\begin{array}{ll}
=\inner{c_s^2\rho\div \u,\div \u_n} - \inner{\rho\opd \u,(\omega+i\Db+i\angvel\times) \u_n}\\
+\inner{\div\u,\nabla p\cdot \u_n}+\inner{\nabla p\cdot\u,\div \u_n}+\inner{(\hess(p)-\rho\hess(\phi)) \u,\u_n} 
- i\omega\inner{\gamma\rho\u,\u_n}
\end{array}
\right\}
\end{aligned}
\\
&\hspace{-2mm}\begin{aligned}\label{eq:limsmall}
\left.\begin{array}{ll}
+\,\inner{c_s^2\rho\div (p_n\u-\u),\div \u_n}\\
-\inner{\rho((\omega+i\angvel\times) (p_n\u-\u)+\Db p_n\u-\conv\u),(\omega+i\Db+i\angvel\times) \u_n}\\
+\,\inner{\div(p_n\u-\u),\nabla p\cdot \u_n}+\inner{\nabla p\cdot(p_n\u-\u),\div \u_n}\\
+\,\inner{(\hess(p)-\rho\hess(\phi)) (p_n\u-\u),\u_n} 
-\,i\omega\inner{\gamma\rho(p_n\u-\u),\u_n}.
\end{array}
\right\}
\end{aligned}
\end{align}
It holds that $\lim_{n\in\IN''}\eqref{eq:lima}=a(\u,\u')$.
Further we estimate that $|\eqref{eq:limsmall}|\lesssim \dn(\u,p_n\u)$ and hence $\lim_{n\in\IN''}\eqref{eq:limsmall}=0$ due to \cref{lem:dnpn}.
Thus altogether we obtain that $\lim_{n\in\IN''} \|A_np_n\u-p_nA\u\|_{\Xn}=0$ and hence $\lim_{n\to\infty} \|A_np_n\u-p_nA\u\|_{\Xn}=0$, which finishes the proof.
\end{proof}

\section{Convergence analysis}\label{sec:convergence}

\subsection{Weak right T-coercivity}

% \marginpar{include: get rid of ``$\dom$ convex''? not necessary for stars, but maybe of general interest}

First we recall how the well posedness of the continuous problem \eqref{eq:Galbrun} is established \cite{HallaHohage:21}.
That is the injectivity of \eqref{eq:Galbrun} \cite[Lemma~3.7]{HallaHohage:21} which follows in a straightforward fashion combined with the weak T-coercivity of \eqref{eq:Galbrun}.
Actually for the latter we use in the current article \emph{right} T-coercivity instead of \emph{left} T-coercivity as in \cite{HallaHohage:21}, and we also choose a slightly different construction of $T$ compared to \cite{HallaHohage:21}.
The reasons for these changes are to be aligned with the forthcoming discrete analysis.
To construct $T$ we first derive a toplogical decomposition of $\IX$.
Thus for $\u\in H_0(\div)$ we seek a solution $\vsi\in H^2$ to
\begin{subequations}\label{eq:def-v-ohne-K}
\begin{align}
(\div+\q\cdot)\nabla\vsi &=(\div+\q\cdot)\u &\text{in }\dom,\\
\nv\cdot\nabla \vsi &= 0 &\text{on }\partial\dom.
\end{align}
\end{subequations}
Note at this point that we only demand $\u\in H_0(\div)\subset\IX$ and the reason why we emphasize this is that the discrete spaces satisfy $\Xn\subset H_0(\div)$, but $\Xn\not\subset\IX$ (in general).
To start with we consider \eqref{eq:def-v-ohne-K} as variational problem in $H^1$.
If a solution $\vsi\in H^1$ to \eqref{eq:def-v-ohne-K} exists, then it follows with convenient regularity theory \cite[Theorem~2.17]{ABDG:98} that the map $\u\to\vsi$ is in $L(H_0(\div),H^2)$.
Although the sesquilinear form associated to the left hand-side of \eqref{eq:def-v-ohne-K} is only weakly coercive and we cannot guarantee the injectivity of the associated operator.
As a remedy we consider the problem on $H^2_*$ and introduce in addition to the low order perturbation $\q\cdot$ another perturbation through an operator $\Kuni$.
In addition we replace $\q\cdot$ by $\PLz\q\cdot$ (with $\PLz$ being the orthogonal projection from $L^2$ to $L^2_0$) to enable a suitable perturbation analysis.
Thus let $\Htbc:=\{\phi\in H^2_*\colon \nv\cdot\nabla \phi = 0\}$ and 
\begin{align*}
\Kuni:=\sum_{l=1}^L \psi_l \spl \div \cdot,\div\nabla\phi_l \spr,
\end{align*}
whereat $L\in\IN_0$ is the dimension of the kernel space of $(\div+\PLz\q\cdot)\nabla\in L(\Htbc,L^2_0)$, $\phi_l\in\Htbc$, $l=1,\dots,L$ is an orthonormal basis with respect to the $\Htbc$-equivalent inner product $\spl \div \cdot,\div\cdot \spr$ of the kernel space, and $\psi_l\in L^2_0$, $l=1,\dots,L$ is an orthonormal basis of the $L^2_0$-orthogonal complement of $(\div+\PLz\q\cdot)\nabla\Htbc$.
Thence we consider the equation
\begin{subequations}\label{eq:def-vi}
\begin{align}
(\div+\PLz\q\cdot+\Kuni)\nabla\vsi &=(\div+\PLz\q\cdot+\Kuni)\u &\text{in }\dom,\\
\nv\cdot\nabla\vsi &= 0 &\text{on }\partial\dom,
\end{align}
\end{subequations}
instead of \eqref{eq:def-v-ohne-K}.
Thus for $\u\in\IX$ we set
\begin{align}\label{eq:Def-PV}
\vi:=\PV\u:=\nabla\vsi,
\end{align}
$\w:=\u-\vi$ and
\begin{align}\label{eq:Def-T}
T:=\vi-\wi.
\end{align}
It follows from its construction that $T\in L(\IX)$ and $TT=\Id_\IX$.
It holds that $(\div+\q\cdot)\wi=(\Id_{L^2}-\PLz)(\q\cdot\wi)-M\wi$.
Since the operators $\Id_{L^2}-\PLz$ and $M$ are compact the proof of \cite[Theorem~]{HallaHohage:21} needs only to be adapted slightly to obtain that $A$ is weakly right $T$-coercive.
We do not give more details at this point, because the proof will be contained in the proof of \cref{thm:dis-wTc}.

\subsection{Construction and properties of $T_n$}

Next we introduce a discrete variant $T_n$ of $T$ and analyze its properties.

\subsubsection{Definition of $T_n$}

For $\u\in H_0(\div)$ let $\vsii\in H^2_*$ be the solution to
\begin{subequations}\label{eq:def-vii}
\begin{align}
(\div+\PLz\q\cdot+\Kuni)\nabla\vsii &=(\div+\pi_n^l\q\cdot+\Kuni)\u,\\
\nv\cdot\nabla\vsii &= 0.
\end{align}
\end{subequations}
Note that in comparison to \eqref{eq:def-vi} in the right hand-side of \eqref{eq:def-vii} the term $\PLz\q\cdot$ is changed to $\pi_n^l\q\cdot$ and therefore we use a new symbol $\vsii$ do denote the solution to \eqref{eq:def-vii}.
Subsequently let
\begin{align}\label{eq:Def-tPV-PV}
\vii:=\tPV\u:=\nabla\vsii \quad\text{and}\quad
\v_n:=\PVn \u:=\pi_n^d \nabla\vsii.
\end{align}
For $\u_n\in\Xn$ let $\w_n:=\u_n-\v_n$ and set
\begin{align}\label{eq:Def-Tn}
T_n\u_n:=\v_n-\w_n.
\end{align}

\subsubsection{Boundedness of $T_n$}

\begin{lemma}\label{lem:PVnbounded}
There exists a constant $C>0$ such that $\|\PVn\|_{L(\Xn)}\leq C$ for all $n\in\IN$.
\end{lemma}
\begin{proof}
Let $\u_n\in\Xn$ be given and $\vsii$ be the solution to \eqref{eq:def-vii}.
First we note that $\|\vsii\|_{H^2}\lesssim \|\u_n\|_{\Xn}$.
Since $\nabla\vsii\in \bH^1$ the function $\pi_n^d\nabla\vsii$ is well defined and the uniform boundedness of $\PVn$ follows from $\div\pi_n^d\nabla\vsii=\pi_n^l\div\nabla\vsii$, the point-wise convergence of $\pi_n^d, \pi_n^l$, \eqref{eq:EstApr}, \eqref{eq:EstFace} and $\bjump{\nabla\vsii}=0$.
\end{proof}
\begin{lemma}\label{lem:Tnbounded}
There exists a constant $C>0$ such that $\|T_n\|_{L(\Xn)}\leq C$ for all $n\in\IN$.
\end{lemma}
\begin{proof}
Follows from \cref{lem:PVnbounded} and the definition \eqref{eq:Def-Tn} of $T_n$.
\end{proof}

\subsubsection{Stability of $T_n$}

\begin{lemma}\label{lem:PVnQuasiProj}
Let $\On:=\PVn\PVn-\PVn$.
It holds that $\lim_{n\to\infty}\|\On\|_{L(\Xn)}=0$.
\end{lemma}
\begin{proof}
Let $\u_n\in\Xn$ and $\vsii_1$ be the solution to \eqref{eq:def-vii}.
Thence $\PVn\u_n=\pi_n^d\nabla\vsii_1$.
Let $\vsii_2$ be the solution to \eqref{eq:def-vii} with $\u_n$ being replaced by $\pi_n^d\nabla \vsii_1$ in the right hand-side.
We compute
\begin{align}\label{eq:Dwn}
\begin{aligned}
(\div+\PLz\q\cdot+\Kuni)\nabla \vsii_2
&= (\div+\pi_n^l\q\cdot+\Kuni)\pi_n^d\nabla \vsii_1\\
&= \pi_n^l(\div+\PLz\q\cdot+\Kuni)\nabla \vsii_1 + \Kuni\pi_n^d\nabla \vsii_1-\pi_n^l\Kuni\nabla \vsii_1\\
&= \pi_n^l(\div+\pi_n^l\q\cdot+\Kuni)\u_n + \Kuni(\pi_n^d-\Id_{\IX})\nabla \vsii_1+ (\Id_{L^2_0}-\pi_n^l)\Kuni\nabla \vsii_1\\
&= (\div+\pi_n^l\q\cdot+\Kuni)\u_n + \tOn\u_n
\end{aligned}
\end{align}
with $\tOn\u_n:=\Kuni(\pi_n^d-\Id_{\IX})\nabla \vsii_1 + (\Id_{L^2_0}-\pi_n^l)\Kuni\nabla \vsii_1
 + (\pi_n^l-\Id_{L^2_0})\Kuni\u_n$.
Since $\Kuni$ is a compact operator which maps into $L^2_0$ and $\Id_{L^2_0}-\pi_n^l$ converges point-wise to zero, it follows that $(\Id_{L^2_0}-\pi_n^l)\Kuni$ tends to zero in operator norm.
For the remaining first term in $\tOn$ we estimate
\begin{align*}
\|\Kuni(\pi_n^d-\Id_{\IX})\tPV\|_{L(\Xn,L^2_0)}
\lesssim \|\div(\pi_n^d-\Id_{\IX})\tPV\|_{L(\Xn,L^2_0)}
= \|(\pi_n^l-\Id_{L^2_0})\div\tPV\|_{L(\Xn,L^2_0)}.
\end{align*}
We compute
\begin{align*}
\div \nabla\vsii&=(\div+\pi_n^l\q\cdot+\Kuni)\bu_n-(\PLz\q\cdot+\Kuni)\nabla\vsii,\\
\pi_n^l\div \nabla\vsii&=(\div+\pi_n^l\q\cdot+\Kuni)\bu_n
-\pi_n^l(\q\cdot+\Kuni)\nabla\vsii+(\pi_n^l-\Id)\Kuni\bu_n.
\end{align*}
Hence
\begin{align*}
\|\Kuni(\pi_n^d-\Id_{\IX})\tPV\|_{L(\Xn,L^2_0)}
\lesssim \|(\PLz-\pi_n^l)(\q\cdot+\Kuni)\tPV\|_{L(H_0(\div),L^2_0)}
+\|(\pi_n^l-\Id)\Kuni\|_{L(H_0(\div),L^2_0)}.
\end{align*}
The former right hand-side tends to zero due to the previously used arguments
and
\begin{align*}
\|(\PLz-\pi_n^l)(\q\cdot\tPV)\|_{L(H_0(\div),L^2_0)}
\lesssim \hmax \|\q\cdot\tPV\|_{L(H_0(\div),H^1)}
\lesssim \hmax.
\end{align*}
Now it follows from $(\div+\PLz\q\cdot+\Kuni)\nabla (\vsii_2-\vsii_1)=\tOn\u_n$ that
\begin{align*}
\|(\PVn\PVn-\PVn)\bu_n\|_{\Xn}
\lesssim \|\nabla\vsii_2-\nabla\vsii_1\|_{\bH^1}
\lesssim \|\tOn\|_{L(\Xn,L^2_0)} \|\bu_n\|_{\Xn}
\end{align*}
and the claim is proven.
\end{proof}
\begin{lemma}
There exist constants $n_0,C>0$ such that $T_n$ is invertible and $\|T_n^{-1}\|_{L(\Xn)}\leq C$ for all $n>n_0$.
\end{lemma}
\begin{proof}
It follows from \cref{lem:PVnQuasiProj} and the definition of $T_n$ that $T_nT_n=4\PVn\PVn-4\PVn+\Id_{\Xn}=\Id_{\Xn}+4\On$.
Since $\|\On\|_{L(\Xn)}$ tends to zero there exists an index $n_0>0$ such that $\|\On\|_{L(\Xn)}<1/8$ for all $n>n_0$ and hence $\|(\Id_{\Xn}+4\On)^{-1}\|_{L(\Xn)}\leq 2$ for $n>n_0$.
It follows that $T_n^{-1}=(\Id_{\Xn}+4\On)^{-1}T_n$ and thus $\|T_n^{-1}\|_{L(\Xn)}\leq 2\|T_n\|_{L(\Xn)}$, $n>n_0$.
The claim follows now from \cref{lem:Tnbounded}.
\end{proof}

\subsubsection{Asymptotic consistency of $T_n$}

\begin{lemma}\label{lem:PVnAsympCons}
It holds that $\lim_{n\to\infty} \|(\PVn p_n-p_n \PV)\u\|_{\Xn}=0$ for each $\u\in\IX$.
\end{lemma}
\begin{proof}
Recall that $\PV, \tPV$ and $\PVn$ are defined in \eqref{eq:Def-PV} and \eqref{eq:Def-tPV-PV} respectively.
We estimate
\begin{align*}
\|(\PVn p_n-p_n \PV)&\u\|_{\Xn}\\
&\leq \dn(\PV\u,p_n \PV\u) + \dn(\PV\u,\PVn p_n \u)\\
&= \dn(\PV\u,p_n \PV\u) + \dn(\PV\u,\pi_n^d \tPV p_n \u)\\
&\leq \dn(\PV\u,p_n \PV\u) + \dn(\PV\u,\pi_n^d \tPV \u) + \|\pi_n^d \tPV (\u-p_n\u)\|_{\Xn}\\
&\lesssim \dn(\PV\u,p_n \PV\u)
+ \dn(P_V\u,\pi_n^d P_V \u)
+ \|\PV\u-\tPV \u\|_{\IX}
+ \|\u-p_n\u\|_{H(\div)}\\
&\lesssim \dn(\PV\u,p_n \PV\u)
+ \dn(P_V\u,\pi_n^d P_V \u)
+ \|(\PLz-\pi_n^l) (\q\cdot\u)\|_{L^2}
+ \dn(\u,p_n\u).
\end{align*}
The claim follows now from \cref{lem:dnpn,lem:dnpi} and the point-wise convergence of $\pi_n^l$.
\end{proof}
\begin{lemma}\label{lem:Tnasympcons}
It holds that $\lim_{n\to\infty} \|(T_n p_n-p_n T)\u\|_{L(\IX,\Xn)}=0$ for each $\u\in\IX$.
\end{lemma}
\begin{proof}
Follows from \cref{lem:PVnAsympCons} and the definitions \eqref{eq:Def-T}, \eqref{eq:Def-Tn} of $T, T_n$.
\end{proof}

\subsection{Discrete weak $T_n$-coercivity}

Let us collect some properties, to emphasize which are the essential ingredients for the proof of the forthcoming \cref{thm:dis-wTc}.
We recall the identity \cite[\bb3.7)]{HallaHohage:21}
\begin{align}\label{eq:Grisvard}
\spl c_s^2\rho \div\v,\div\v \spr
&=|\v|_{\bH^1_{c_s^2\rho}}^2  + \spl \Keta\v,\Keta\v \spr_{\bV}
\end{align}
with a compact operator $\Keta\in L(\bV)$ for all $\v\in\bV$, whereat
\begin{align*}
\bV:=\{\nabla v\colon v\in\Htbc\}, \quad
\|\cdot\|_{\bV}:=|\cdot|_{\bH^1_{c_s^2\rho}}.
\end{align*}
Note that actually \eqref{eq:Grisvard} is in \cite{HallaHohage:21} formulated only for a subspace of $\bV$.
However, the proof to extend this result to $\bV$ requires no changes at all.
The purpose of the next lemma is to work out an estimate in weighted norms, which is essential to obtain the robustness with respect to $\rhol/\rhou$ and $\csl/\csu$.
\begin{lemma}\label{lem:Est-pin}
It holds that
\begin{align*}
\| \rho^{\sfrac12} \Db \pi_n^d\v \|_{\bL^2}^2
\leq (\Cpis)^2 (1+\hmax^2\Cpih) \Mach^2 |\v|_{\bH^1_{c_s^2\rho}}^2
\end{align*}
with constants $\Cpih>0$,
\begin{align*}
(\Cpis)^2:=2\big((\Cface\Chratio\Ctr)^2+\sup_{n\in\IN}\sup_{\tet\in\calT_n}\|\pi_n^d\|_{L(\bH^1_*(\tet))}^2\big), \qquad
\|\cdot\|_{\bH^1_*(\tet)}:=|\cdot|_{\bH^1(\tet)}
\end{align*}
for all $\v\in\Hz$, $n\in\IN$.
\end{lemma}
\begin{proof}
For each $\tet\in\calT_n$ we estimate that
\begin{align*}
\|\rho^{\sfrac12}\conv \pi_n^d\v\|_{\bL^2(\tet)}^2
&\leq \|c_s^{-1}\bflow\|_{\bL^\infty}^2 \csu_\tet^2\rhou_\tet |\pi_n^d\v|_{\bH^1(\tet)}^2\\
&\leq \|c_s^{-1}\bflow\|_{\bL^\infty}^2 \csu_\tet^2\rhou_\tet \|\pi_n^d\|_{L(\bH^1_*(\tet))}^2 \|\v|_{\bH^1(\tet)}^2\\
&\leq \|c_s^{-1}\bflow\|_{\bL^\infty}^2 \|\pi_n^d\|_{L(\bH^1_*(\tet))}^2
\Big(1+\hmax^2 \frac{1}{\csl^2\rhol} (\CL_{c_s\rho^{\sfrac12}})^2\Big)^2
|\v|_{\bH^1_{c_s^2\rho}(\tet)}^2.
\end{align*}
We further compute that
\begin{align*}
\|\rho^{\sfrac12}\LO\pi_n^d\v\|_{\bL^2(\tet)}
=\|\rho^{\sfrac12}\sum_{\face\in\calF_\tet}\LO^\face\pi_n^d\v\|_{\bL^2(\tet)}
&\leq \Ctr \rhou_\tet \|\h^{\sfrac{-1}{2}}\bjump{\pi_n^d\v}\|_{\bL^2(\partial\tet)}\\
&=\Ctr \rhou_\tet \|\h^{\sfrac{-1}{2}}\bjump{\pi_n^d\v-\v}\|_{\bL^2(\partial\tet)}.
\end{align*}
Hence we estimate by means of \eqref{eq:EstFace} that
\begin{align*}
\Ctr^2 \sum_{\tet\in\calT_n} \rhou_\tet \|\h^{\sfrac{-1}{2}}\bjump{\pi_n^d\v-\v}\|_{\bL^2(\partial\tet)}^2
&\leq \Ctr^2 \sum_{\tet\in\calT_n} \rhou_\tet \sum_{\face\in\calF_\tet} \|\h^{\sfrac{-1}{2}}\bjump{\pi_n^d\v-\v}\|_{\bL^2(\face)}^2\\
&\leq \Ctr^2 \sum_{\tet\in\calT_n} \rhou_\tet \sum_{\face\in\calF_\tet} \Big(\frac12 \sum_{j=1}^2 \|\h^{\sfrac{-1}{2}}(\nv\cdot\vel)((\pi_n^d\v)_j-\v)\|_{\bL^2(\face)}\Big)^2\\
&\leq \frac{\Ctr^2}{2} \sum_{\tet\in\calT_n} \rhou_\tet \sum_{\face\in\calF_\tet}  \sum_{j=1}^2 \|\h^{\sfrac{-1}{2}}(\nv\cdot\vel)((\pi_n^d\v)_j-\v)\|_{\bL^2(\face)}^2\\
&\leq \Mach^2 \Ctr^2 \sum_{\tet\in\calT_n} \csu_\tet^2 \rhou_{\dom_\tet} \|\h^{\sfrac{-1}{2}}((\pi_n^d\v)|_{\tet}-\v)\|_{\bL^2(\partial\tet)}^2\\
&\leq\Mach^2\Cface^2\Chratio^2\Ctr^2 \sum_{\tet\in\calT_n} \csu_\tet^2 \rhou_{\dom_\tet} |\v|_{\bH^1(\tet)}^2\\
&\leq\Mach^2\Cface^2\Chratio^2\Ctr^2
\Big(1+(C\hmax)^2 \frac{1}{\csl^2\rhol} (\CL_{c_s\rho^{\sfrac12}})^2\Big)^2 |\v|_{\bH^1_{c_s^2\rho}}^2
\end{align*}
with a constant $C>0$ that only depends on $\Chratio$.
Thus the claim follows.
\end{proof}
It further follows from \eqref{eq:Dwn} that for $\u_n\in\Xn$
\begin{align}\label{eq:divq-wn}
(\div+\pi_n^l\q\cdot)\w_n=-\Kuni\w_n-\tOn\u_n,
\end{align}
whereat we recall that $\Kuni\in L(H_0(\div),L^2_0)$ is compact and $\|\tOn\|_{L(\Xn,L^2_0)}$ tends to zero.
The next lemma shows that operators such as $\Kuni$ lead indeed to compact sequences of operators in the sense of discrete approximation schemes.
\begin{lemma}\label{lem:Kncompact}
Let $K_n^{E\PV}, K_n^{\Keta}, K_n^{\Kuni} \in L(\Xn)$ be defined by
\begin{align*}
\spl K_n^{E\PV}\u_n,\u_n'\spr_{\Xn}&:=\spl \PVn\u_n,\PVn\u_n'\spr_{\bL^2},\quad
\spl K_n^{\mean}\u_n,\u_n'\spr_{\Xn}:=\spl \mean(\q\cdot\w_n),\mean(\q\cdot\w_n')\spr_{L^2},\\
\spl K_n^{\Keta}\u_n,\u_n'\spr_{\Xn}&:=\spl \Keta\tPV\u_n,\Keta\tPV\u_n'\spr_{\bV},\quad
\spl K_n^{\Kuni}\u_n,\u_n'\spr_{\Xn}:=\spl \Kuni\u_n,\Kuni\u_n'\spr_{L^2},
\end{align*}
for all $\u_n,\u_n'\in\Xn$.
Then $(K_n^{E\PV})_{n\in\IN}$, $(K_n^{\Keta})_{n\in\IN}$, $(K_n^{\Kuni})_{n\in\IN}$ are compact in the sense of discrete approximation schemes.
\end{lemma}
\begin{proof}
Let $(\u_n)_{n\in\IN}$, $\u_n\in\Xn$ be a given bounded sequence $\|\u_n\|_{\Xn}\leq1$ for each $n\in\IN$.
Let an arbitrary subsequence $\IN'\subset\IN$ be given.
Recall the compact Sobolev embedding $E_{\Hz,\bL^2}\in L(\Hz,\bL^2)$ and that $\PVn=\pi_n^d\tPV$, $\tPV\in L(H_0(\div),\Hz)$, $n\in\IN$ are uniformly bounded.
Thence there exists $\z\in\bL^2$ and a subsequence $\IN''\subset\IN'$ such that $\lim_{n\in\IN''}\|\z-\tPV\u_n\|_{\bL^2}=0$.
It further follows that
\begin{align*}
\|\z-\PVn\u_n\|_{\bL^2} = \|\z-\pi_n^d\tPV\u_n\|_{\bL^2}
&\leq \|\z-\tPV\u_n\|_{\bL^2}+\|(1-\pi_n^d)\tPV\u_n\|_{\bL^2}\\
&\lesssim \|\z-\tPV\u_n\|_{\bL^2} + \hmax \|\tPV\u_n\|_{\bH^1}
\xrightarrow{n\in\IN''}0.
\end{align*}
We want to show that $\lim_{n\in\IN'''}\|p_n\PV^*\z-K_n^{E\PV}\u_n\|_{\Xn}=0$ for a subsequence $\IN'''\subset\IN''$.
To this end let $\u_n'\in\Xn$, $\|\u_n'\|_{\Xn}=1$, $n\in\IN''$ be such that $\|p_n\PV^*\z-K_n^{E\PV}\u_n\|_{\Xn}\leq |\spl p_n\PV^*\z-K_n^{E\PV}\u_n,\u_n'\spr_{\Xn}|+1/n$.
By means of \cref{lem:weakcompactness} we choose $\IN'''\subset\IN''$ and $\u\in\IX$ such that $\u_n'$ converges weakly to $\u$ in the sense of \cref{lem:weakcompactness}.
We compute
\begin{align*}
\spl p_n \PV^*\z,\u_n'\spr_{\Xn}
&=\spl \div \PV^*\z,\div\u_n'\spr+\spl \PV^*\z,\u_n'\spr+\spl \conv \PV^*\z,\Db\u_n'\spr\\
&\xrightarrow{n\in\IN'''} \spl \div \PV^*\z,\div\u\spr+\spl \PV^*\z,\u\spr+\spl \conv \PV^*\z,\conv\u\spr
=\spl \PV^*\z,\u\spr_{\IX}
=\spl \z,\PV\u\spr
\end{align*}
and
\begin{align*}
\spl K_n^{E\PV}\u_n,\u_n'\spr_{\Xn}
&=\spl \PVn\u_n,\PVn\u_n'\spr_{\bL^2}\\
&=\spl \PVn\u_n-\z,\PVn\u_n'\spr_{\bL^2}+\spl \z,\pi_n^d\tPV\u_n'\spr_{\bL^2}\\
&=\spl \PVn\u_n-\z,\PVn\u_n'\spr_{\bL^2}+\spl \z,\pi_n^d\PV\u_n'\spr_{\bL^2}
+\spl \z,\pi_n^d(\tPV-\PV)\u_n'\spr_{\bL^2}\\
&=\spl \PVn\u_n-\z,\PVn\u_n'\spr_{\bL^2}
+\spl \z,\PV\u_n'\spr_{\bL^2}
+\spl \z,(\pi_n^d-1)\PV\u_n'\spr_{\bL^2}\\
&+\spl \z,\pi_n^d(\tPV-\PV)\u_n'\spr_{\bL^2}\\
&=\spl \PVn\u_n-\z,\PVn\u_n'\spr_{\bL^2}
+\spl \z,\PV\u_n'\spr_{\bL^2}
+\spl \z,(\pi_n^d-1)\PV\u_n'\spr_{\bL^2}\\
&+\spl \z,(\pi_n^d-1)(\tPV-\PV)\u_n'\spr_{\bL^2}
+\spl \z,(\tPV-\PV)\u_n'\spr_{\bL^2}.
\end{align*}
As previously we estimate
\begin{align*}
\|(1-\pi_n^d)\PV\u_n'\|_{\bL^2}+\|(1-\pi_n^d)(\tPV-\PV)\u_n'\|_{\bL^2}
\lesssim \hmax (\|\PV\u_n'\|_{\bH^1}+\|\tPV\u_n'\|_{\bH^1}) \lesssim \hmax.
\end{align*}
In addition we can write $(\tPV-\PV)\u_n'=S(\pi_n^l-\PLz)(\q\cdot\u_n')$ with $S:=\nabla\big((\div+\PLz\q\cdot+\Kuni)\nabla\big)^{-1} \in L(L^2_0,\bL^2)$ and hence
\begin{align*}
\spl \z,(\tPV-\PV)\u_n'\spr_{\bL^2}
&=\spl \z,S(\pi_n^l-\PLz)(\q\cdot\u_n')\spr_{\bL^2}
=\spl (\pi_n^l-\PLz)S^*\z,\PLz(\q\cdot\u_n')\spr_{L^2_0}
\xrightarrow{n\in\IN'''}0,
\end{align*}
whereat we used that $\pi_n^l$ is indeed an orthogonal projection which converges point-wise.
Since
\begin{align*}
\spl \z,\PV\u_n'\spr_{\bL^2}
=\spl \PV^*\z,\u_n'\spr_{H_0(\div)}
\xrightarrow{n\in\IN'''}\spl \PV^*\z,\u\spr_{H_0(\div)}
=\spl \z,\PV\u\spr_{\bL^2}
\end{align*}
the claim for $(K_n^{E\PV})_{n\in\IN}$ follows.
The proofs for $(K_n^{\mean})_{n\in\IN}$, $(K_n^{\Keta})_{n\in\IN}$ and $(K_n^{\Kuni})_{n\in\IN}$ can be derived by the very same technique.
\end{proof}
In order to formulate Theorem~\ref{thm:dis-wTc} we introduce some additional quantities.
Let $\lambda_-(\m) \in L^\infty$ be the smallest eigenvalue of a positive definite matrix and $\m:=-\rho^{-1}\hess(p)+\hess(\phi)$.
Further let
\begin{align*}%\label{eq:theta}
\Cm&:=\max\Big\{0, \sup_{x\in\dom} \frac{-\lambda_-(\m(x))}{\gamma(x)}\Big\}
\qquad\text{and}\qquad
\theta:=\arctan(\Cm/|\omega|)\in[0,\pi/2), \quad \text{for }\omega\neq0.
\end{align*}
\begin{theorem}\label{thm:dis-wTc}
If $\Mach^2<\frac{1}{(\Cpis)^2} \frac{1}{1+\Cm^2/|\omega|^2}$, then
$A_nT_n=B_n+K_n$ with $(B_n\in L(\Xn))_{n\in\IN}$ being uniformly bounded and stable, $(K_n\in L(\Xn))_{n\in\IN}$ being compact, and there exists a bijective operator $B\in L(\IX)$ such that $\lim_{n\to\infty} \|B_n p_n \u - p_n B \u\|_{\Xn}=0$ for each $\u\in\IX$.
\end{theorem}
\begin{proof}
\emph{1.\ step: definition of $B_n$ and $K_n$.}\quad
Let $\Keta\in L(\bV)$ be the compact operator from \eqref{eq:Grisvard}.
Let
\begin{subequations}\label{eq:tBn}
\begin{align}
% v_n, v_n'
\nonumber
\spl \tilde B_n\u_n,&\u'_n \spr_{\Xn}:=\\
&\label{eq:tBn-a}
\inner{c_s^2\rho\div \v_n,\div \v_n'}
-\inner{\rho i\Db\v_n,i\Db \v_n'}
+\inner{c_s^2\rho \pi_n^l(\q\cdot \w_n),\pi_n^l(\q\cdot \w_n')}\\
% v_n, w_n', w_n, v_n'
&
-\inner{\rho i\Db\v_n,\opdn \w_n'}
+\inner{\rho \opdn \w_n,i\Db\v_n'}\\
% w_n, w_n'
&
+\inner{\rho\opdn \w_n,\opdn \w_n'}
+\inner{\rho(i\gamma+\m) \w_n,\w_n'}\\
&\label{eq:tBn-stab}
+\inner{\v_n,\v_n'}
+\CKeta \spl \Keta\tPV\u_n,\Keta\tPV\u_n' \spr_{\bV}
+\inner{\Kuni\w_n,\Kuni \w_n'}
+\inner{\tOn\u_n,\tOn \u_n'}
\end{align}
\end{subequations}
and
\begin{subequations}\label{eq:tKn}
\begin{align}
\nonumber
\spl \tilde K_n&\u_n,\u'_n \spr_{\Xn}:=\\
&\label{eq:tKn-stabi}
\CtKn (\spl\v_n,\v_n'\spr
+\spl \Keta\tPV\u_n,\Keta\tPV\u_n' \spr_{\bV}
+\spl \tOn\u_n,\tOn\u_n' \spr\\
&\label{eq:tKn-stabii}
+\spl\Kuni\w_n,\Kuni \w_n'\spr
+\spl \mean(\q\cdot\w_n), \mean(\q\cdot\w_n') \spr)\\
% v_n, v_n'
\nonumber\\
&\label{eq:tKn-b}
+\inner{c_s^2\rho\q\cdot\v_n,\div\v_n'}
+\inner{c_s^2\rho\div\v_n,\q\cdot\v_n'}
-\inner{\rho (\omega+i\angvel\times) \v_n, (\omega+i\angvel\times) \v_n'}\\
&\label{eq:tKn-c}
-\inner{\rho (\omega+i\angvel\times) \v_n, i\Db \v_n'}
-\inner{\rho i\Db \v_n, (\omega+i\angvel\times) \v_n'}
-i\omega\inner{\gamma\rho\v_n,\v_n'}\\
&\label{eq:tKn-d}
-\inner{\rho\m \v_n,\v_n'}\\
% v_n, w_n'
\nonumber\\
&\label{eq:tKn-e}
-\inner{\rho\m \v_n,\w_n'}
-i\omega\inner{\gamma\rho\v_n,\w_n'}
-\inner{c_s^2\rho \pi_n^l(\q\cdot \v_n),\pi_n^l(\q\cdot \w_n')}\\
&\label{eq:tKn-f}
-\inner{\rho (\omega+i\angvel\times) \v_n, \opdn \w_n'}\\
&\label{eq:tKn-g}
-\inner{c_s^2\rho (\div+\pi_n^l\q\cdot)\v_n,\Kuni\w_n'+\tOn\u_n'}
+\spl c_s^2\rho (\Id-\pi_n^l)(\q\cdot\v_n) , \div\w_n' \spr\\
% w_n, v_n'
\nonumber\\
&\label{eq:tKn-h}
+\inner{\rho\m \w_n,\v_n'}
+i\omega\inner{\gamma\rho\w_n,\v_n'}
+\inner{c_s^2\rho \pi_n^l(\q\cdot \w_n), \pi_n^l(\q\cdot \v_n')}\\
&\label{eq:tKn-i}
+\inner{\rho \opdn \w_n, (\omega+i\angvel\times) \v_n'}\\
&\label{eq:tKn-j}
+\inner{c_s^2\rho (\Kuni\w_n+\tOn\u_n), (\div+\pi_n^l\q\cdot)\v_n'}
-\spl c_s^2\rho \div\w_n, (\Id-\pi_n^l)(\q\cdot\v_n') \spr\\
% w_n, w_n'
\nonumber\\
&\label{eq:tKn-k}
-\spl c_s^2\rho (\Id-\mean-\pi_n^l)(\q\cdot\w_n) , \div\w_n' \spr
-\spl c_s^2\rho \mean(\q\cdot\w_n) , \div\w_n' \spr\\
&\label{eq:tKn-l}
-\spl c_s^2\rho \div\w_n, (\Id-\mean-\pi_n^l)(\q\cdot\w_n') \spr
-\spl c_s^2\rho \div\w_n, \mean(\q\cdot\w_n') \spr\\
&\label{eq:tKn-m}
-\spl c_s^2\rho (\Kuni\w_n+\tOn\u_n),\Kuni\w_n'+\tOn\u_n' \spr 
\end{align}
\end{subequations}
and
\begin{subequations}\label{eq:Kn}
\begin{align}
\nonumber
\spl K_n&\u_n,\u'_n \spr_{\Xn}:=\\
&\label{eq:Kn-a}
-\CtKn \spl\v_n,\v_n'\spr
-(\CKeta+\CtKn) \spl \Keta\tPV\u_n,\Keta\tPV\u_n' \spr_{\bV}\\
&\label{eq:Kn-b}
-(1+\CtKn) \spl\Kuni\w_n,\Kuni \w_n'\spr
-\CtKn \spl \mean(\q\cdot\w_n), \mean(\q\cdot\w_n') \spr
-(1+\CtKn)\inner{\tOn\u_n,\tOn \u_n'}
\end{align}
\end{subequations}
for all $\u_n,\u_n'\in\Xn$, whereat the constansts $\CKeta,\CtKn>0$ will be specified lateron.
Thence with $B_n:=\tilde B_n+\tilde K_n$ it holds that $A_nT_n=B_n+K_n$.
We discuss the details of this decomposition in the following.
First note that $\spl A_nT_n\u_n,\u_n' \spr_{\Xn}=\sesn(T_n\u_n,\u_n')=\sesn(\v_n-\w_n,\v_n'+\w_n')$.
Second note that the operators $\tilde K_n$ and $K_n$ contain only terms which are intuitively compact.
Although to avoid lengthy and technical proofs that all these terms yield compact sequences of operators in the sense of discrete approximation schemes, we will choose a sufficiently large constant $\CtKn$ such that $\tilde K_n$ becomes small in a suitable sense.
On the other hand the compactness of $(K_n)_{n\in\IN}$ is ensured by \cref{lem:Kncompact} and $\lim_{n\to\infty}\|\tOn\|_{L(\Xn,L^2)}=0$.
Note that we added the line \eqref{eq:tBn-stab} into the definition of $\tilde B_n$ for stability reasons and together with \eqref{eq:tKn-stabi}, \eqref{eq:tKn-stabii} they cancel out with \eqref{eq:Kn-a}, \eqref{eq:Kn-b}.
The term $\inner{c_s^2\rho \pi_n^l(\q\cdot \w_n),\pi_n^l(\q\cdot \w_n')}$ in \eqref{eq:tBn-a} is also included for stability reasons and it is substracted in \eqref{eq:tKn} to keep a zero balance - however, its appearence in \eqref{eq:tKn} is hidden and will be explained shortly.
Also note that in \eqref{eq:tKn} we grouped together all terms of the kinds $(\v_n,\v_n')$, $(\v_n,\w_n')$, $(\w_n,\v_n')$, $(\w_n,\w_n')$ respectively and we included blank lines to emphasize those different blocks.
To get rid of the term $-\spl c_s^2\rho \div\w_n,\div\w_n' \spr$ (which appears naturally in $\sesn(\v_n-\w_n,\v_n'+\w_n')$) we add
\begin{align*}
-\spl c_s^2\rho \div\w_n,\pi_n^l(\q\cdot\w_n') \spr
-\spl c_s^2\rho \pi_n^l(\q\cdot\w_n), \div\w_n' \spr
-\spl c_s^2\rho \pi_n^l(\q\cdot\w_n), \pi_n^l(\q\cdot\w_n') \spr
\end{align*}
to obtain the square
\begin{align*}
-\spl c_s^2\rho (\div+\pi_n^l\q\cdot)\w_n, (\div+\pi_n^l\q\cdot)\w_n' \spr
=-\spl c_s^2\rho (\Kuni\w_n+\tOn\u_n),\Kuni\w_n'+\tOn\u_n' \spr
=\eqref{eq:tKn-m},
\end{align*}
whereat the first equalitiy is due to \eqref{eq:divq-wn}.
Hence we explained the previously stated appearence of $-\spl c_s^2\rho \pi_n^l(\q\cdot\w_n), \pi_n^l(\q\cdot\w_n') \spr$ in \eqref{eq:tKn} which balances to zero with the last term in \eqref{eq:tBn-a}.
It remains to balance the previously added terms $-\spl c_s^2\rho \div\w_n,\pi_n^l(\q\cdot\w_n') \spr
-\spl c_s^2\rho \pi_n^l(\q\cdot\w_n), \div\w_n' \spr$.
With a reverse sign those terms are grouped in \eqref{eq:tKn-k}-\eqref{eq:tKn-l} together with $-\spl c_s^2\rho \div\w_n,\q\cdot\w_n' \spr
-\spl c_s^2\rho \q\cdot\w_n, \div\w_n' \spr$ which appear naturally in $\sesn(\v_n-\w_n,\v_n'+\w_n')$.
In addition the terms
$\mp\spl c_s^2\rho \mean(\q\cdot\w_n) , \div\w_n' \spr$
and
$\mp\spl c_s^2\rho \div\w_n, \mean(\q\cdot\w_n') \spr$
are included in \eqref{eq:tKn-k}, \eqref{eq:tKn-l} respectively for later use.
The terms
$\spl c_s^2\rho \div\v_n,(\div+\q\cdot)\w_n')\spr
-\spl c_s^2\rho (\div+\q\cdot)\w_n),\div\v_n'\spr$
which appear naturally in $\sesn(\v_n-\w_n,\v_n'+\w_n')$ are reformulated in \eqref{eq:tKn-e}-\eqref{eq:tKn-g} and \eqref{eq:tKn-h}-\eqref{eq:tKn-j} the same way as applied previously for \eqref{eq:tKn-k}-\eqref{eq:tKn-m}.
Finally we note that all terms in \eqref{eq:tKn-b}-\eqref{eq:tKn-d} occur naturally in $\sesn(\v_n-\w_n,\v_n'+\w_n')$.
The uniform boundedness of $B_n$, $n\in\IN$ follows straightforwardly.
\\

\emph{2.\ step: coercivity of $\tilde B_n$.}\quad
The commutation $\div\pi_n^d=\pi_n^l\div$ will enable us to adapt \eqref{eq:Grisvard} to $\v_n$ in an apt way.
To this end we compute
\begin{align*}
\div\v_n
&=\div\pi_n^d\nabla\vsii
=\pi_n^l\Delta\vsii
=\pi_n^l\left( -(\PLz\q\cdot+\Kuni)\nabla\vsii+(\div+\pi_n^l\q\cdot+\Kuni)\bu_n \right)\\
&=-(\PLz\q\cdot+\Kuni)\nabla\vsii+(\div+\pi_n^l\q\cdot+\Kuni)\bu_n
+(\Id-\pi_n^l)(\PLz\q\cdot+\Kuni)\nabla\vsii+(\pi_n^l-\Id)M\u_n\\
&=\Delta\vsii
+(\Id-\pi_n^l)(\PLz\q\cdot+\Kuni)\tPV\u_n+(\pi_n^l-\Id)M\u_n\\
&=:\Delta\vsii+\hOn\bu_n.
\end{align*}
By the same technique as used in the proof of \cref{lem:PVnQuasiProj} it follows that $\lim_{n\to\infty}\|\hOn\|_{L(\Xn,L^2_0)}=0$.
Thus with $\spl \cOn\u_n,\u_n' \spr_{\Xn}:=\spl c_s^2\rho \div\v_n,\hOn\u_n' \spr
+\spl c_s^2\rho \hOn\u_n, \div\v_n' \spr
+\spl c_s^2\rho \hOn\u_n, \hOn\u_n' \spr$
it holds that $\lim_{n\to\infty}\|\cOn\|_{\Xn}=0$ and
\begin{align}\label{eq:divdiv-vn}
\spl c_s^2\rho \div\v_n,\div\v_n \spr
=\spl c_s^2\rho \Delta\vsii,\Delta\vsii \spr+\spl \cOn\u_n,\u_n \spr_{\Xn}.
\end{align}
Thus \eqref{eq:divdiv-vn} and \eqref{eq:Grisvard} yield that
\begin{align}\label{eq:divdiv-vnii}
\spl c_s^2\rho \div\v_n,\div\v_n \spr
&=|\nabla\vsii|_{\bH^1_{c_s^2\rho}}^2
+\spl \Keta\tPV\u_n,\tPV\u_n \spr_{\bV}
+\spl \cOn\u_n,\u_n \spr_{\Xn}.
\end{align}
Due to the smallness assumption on the Mach number there exist $\epsilon\in(0,1)$, $\tau\in(0,\pi/2-\theta)$ and $n_0>0$ such that
\begin{align*}
C_{\theta,\tau,\epsilon,n_0}:=1 - (\Cpis)^2 (1+\sup_{n>n_0}\hmax^2\Cpih) \Mach^2 \big(1+\tan^2(\theta+\tau)(1-\epsilon)^{-1}\big)-2\epsilon>0.
\end{align*}
Henceforth we assume that $n>n_0$.
Now we estimate by means of a weighted Young's inequality and the definition of $\theta$ that
\begin{align*}
\frac{1}{\cos(\theta+\tau)}
&\Re\left(e^{-i(\theta+\tau)\sign\omega} \spl\tilde B_n\u_n,\u_n\spr_{\Xn}\right)\\
&\hspace{-10mm}=\|c_s\rho^{\sfrac12}\div\v_n\|_{L^2}^2
-\|\rho^{\sfrac12}\Db\v_n\|_{\bL^2}^2
+\|\v_n\|_{\bL^2}^2
+\CKeta \|\Keta\tPV\bu_n\|_{\bV}^2
+\|\Kuni\w_n\|_{L^2}^2\\
&\hspace{-10mm}
+\|\tOn\w_n\|_{L^2}^2
+\|c_s\rho^{\sfrac12} \pi_n^l(\q\cdot\w_n)\|_{\bL^2}^2
+\|\rho^{\sfrac12}\opdn \w_n\|_{\bL^2}^2
+\spl\rho\m\w_n,\w_n\spr_{\bL^2}\\
&\hspace{-10mm}
+2\tan(\theta+\tau)\sign\omega\Im(\inner{\rho \opdn \w_n,i\Db\v_n})
-|\omega| \tan(\theta+\tau) \|(\gamma\rho)^{\sfrac12}\w_n\|_{\bL^2}^2\\
&\hspace{-10mm}\geq\|c_s\rho^{\sfrac12}\div\v_n\|_{L^2}^2
-\big(1+\tan^2(\theta+\tau)(1-\epsilon)^{-1}\big)\|\rho^{\sfrac12}\Db\v_n\|_{\bL^2}^2
+\|\v_n\|_{\bL^2}^2\\
&\hspace{-10mm}
+\CKeta \|\Keta\tPV\bu_n\|_{\bV}^2
+\|\Kuni\w_n\|_{L^2}^2
+\|\tOn\w_n\|_{L^2}^2
+\|c_s\rho^{\sfrac12} \pi_n^l(\q\cdot\w_n)\|_{\bL^2}^2\\
&\hspace{-10mm}
+\epsilon\|\rho^{\sfrac12}\opdn \w_n\|_{\bL^2}^2
+|\omega| \big(\tan(\theta+\tau)-\tan\theta\big) \|(\gamma\rho)^{\sfrac12}\w_n\|_{\bL^2}^2.
\end{align*}
Thence \Cref{lem:Est-pin} and \eqref{eq:divdiv-vnii} yield that
\begin{align*}
\|c_s\rho^{\sfrac12}&\div\v_n\|_{L^2}^2
-\big(1+\tan^2(\theta+\tau)(1-\epsilon)^{-1}\big)
\|\rho^{\sfrac12}\Db\v_n\|_{\bL^2}^2
+\|\v_n\|_{\bL^2}^2
+\CKeta \|\Keta\tPV\bu_n\|_{\bV}^2\\
&\geq
\epsilon\big(\|c_s\rho^{\sfrac12}\div\v_n\|_{L^2}^2+\|\rho^{\sfrac12}\Db\v_n\|_{\bL^2}^2\big)
+C_{\theta,\tau,\epsilon,n_0}
|\nabla\vsii|_{\bH^1_{c_s^2\rho}}^2
+\|\v_n\|_{\bL^2}^2\\
&+\Big(\CKeta-\frac{1}{4\delta}\Big) \|\Keta\tPV\bu_n\|_{\bV}^2
-\big(\delta\sup_{m\in\IN}\|P_{\tilde V_m}\|_{L(\IX_m,\bV)}^2+\|\cOn\|_{L(\Xn)}\big) \|\u_n\|_{\Xn}^2\\
&\geq \epsilon \min\{\csl^2\rhol,\rhol,1\} \|\v_n\|_{\Xn}^2
+\Big(\CKeta-\frac{1}{4\delta}\Big) \|\Keta\tPV\bu_n\|_{\bV}^2\\
&-\big(\delta\sup_{m\in\IN}\|P_{\tilde V_m}\|_{L(\IX_m,\bV)}^2+\|\cOn\|_{L(\Xn)}\big) \|\u_n\|_{\Xn}^2.
\end{align*}
Further \eqref{eq:divq-wn} yields that
\begin{align*}
4 \big(\|\Kuni\w_n\|_{L^2}^2+\|\tOn\w_n\|_{L^2}^2
+\|c_s\rho^{\sfrac12} \pi_n^l(\q\cdot\w_n)\|_{\bL^2}^2 \big)
\geq \|\div\w_n\|_{L^2}^2
\end{align*}
and hence
\begin{align*}
\|\Kuni\w_n\|_{L^2}^2
&+\|\tOn\w_n\|_{L^2}^2
+\|c_s\rho^{\sfrac12} \pi_n^l(\q\cdot\w_n)\|_{\bL^2}^2
+\epsilon\|\rho^{\sfrac12}\opdn \w_n\|_{\bL^2}^2\\
&+|\omega| \big(\tan(\theta+\tau)-\tan\theta\big) \|(\gamma\rho)^{\sfrac12}\w_n\|_{\bL^2}^2
\gtrsim \|\w_n\|_{\Xn}^2.
\end{align*}
Thus
\begin{align*}
\frac{1}{\cos(\theta+\tau)}
&\Re\left(e^{-i(\theta+\tau)\sign\omega} \spl\tilde B_n\u_n,\u_n\spr_{\Xn}\right)\\
&\hspace{-10mm}\geq \CtB \|\u_n\|_{\Xn}^2
+\Big(\CKeta-\frac{1}{4\delta}\Big) \|\Keta\tPV\bu_n\|_{\bV}^2
-\big(\delta\sup_{m\in\IN}\|P_{\tilde V_m}\|_{L(\IX_m,\bV)}^2+\|\cOn\|_{L(\Xn)}\big) \|\u_n\|_{\Xn}^2
\end{align*}
with a constant $\CtB>0$ independent of $\delta, \CKeta, n>n_0$.
Hence we can choose $\delta>0$ and $n_1>n_0$ such that
\begin{align*}
\frac{1}{\cos(\theta+\tau)}
\Re\left(e^{-i(\theta+\tau)\sign\omega} \spl\tilde B_n\u_n,\u_n\spr_{\Xn}\right)
\geq \frac{\CtB}{2} \|\u_n\|_{\Xn}^2
+\Big(\CKeta-\frac{1}{4\delta}\Big) \|\Keta\tPV\bu_n\|_{\bV}^2
\end{align*}
for all $n>n_1$.
Now we choose $\CKeta>1/(4\delta)$ to obtain that
\begin{align*}
\frac{1}{\cos(\theta+\tau)}
\Re\left(e^{-i(\theta+\tau)\sign\omega} \spl\tilde B_n\u_n,\u_n\spr_{\Xn}\right)
\geq \frac{\CtB}{2} \|\u_n\|_{\Xn}^2
\end{align*}
for all $n>n_1$.
\\

\emph{3.\ step: coercivity of $B_n$.}\quad
To start with we estimate the first term in \eqref{eq:tKn-k} and \eqref{eq:tKn-l} respectively.
To this end we compute that
\begin{align*}
|\spl c_s^2\rho (\Id-\mean-\pi_n^l)(\q\cdot\w_n) , \div\w_n' \spr|
&=|\spl \q\cdot\w_n , (\Id-\mean-\pi_n^l)(c_s^2\rho \div\w_n') \spr|\\
&\leq \|\q\|_{\bL^\infty} \|\w_n\|_{\bL^2} \|(\Id-\mean-\pi_n^l)(c_s^2\rho \div\w_n')\|_{L^2}.
\end{align*}
Thence we apply a discrete commutator technique \cite{Bertoluzza99} and estimate
\begin{align*}
\|(\Id-\mean-\pi_n^l)(c_s^2\rho \div\w_n)\|_{L^2}^2
&=\sum_{\tet\in\calT_n} \|(\Id-\mean-\pi_n^l)(c_s^2\rho \div\w_n)\|_{L^2(\tet)}^2\\
&=\sum_{\tet\in\calT_n} \|(\Id-\mean-\pi_n^l)((c_s^2\rho-c_\tet) \div\w_n)\|_{L^2(\tet)}^2\\
&\leq\sum_{\tet\in\calT_n} \|(c_s^2\rho-c_\tet) \div\w_n\|_{L^2(\tet)}^2\\
&\leq (\CL_{c_s^2\rho})^2 \hmax^2 \sum_{\tet\in\calT_n} \|\div\w_n\|_{L^2(\tet)}^2
= (\CL_{c_s^2\rho})^2 \hmax^2 \|\div\w_n\|_{L^2}^2
\end{align*}
with suitably chosen constants $c_\tet$, $\tet\in\calT_n$.
Let
\begin{align*}
|\u_n|_{Y_n}^2:=
\|\v_n\|_{\bL^2}^2
+\|\Keta\tPV\u_n\|_{\bV}^2
+\|\tOn\u_n\|_{L^2}^2
+\|\Kuni\w_n\|_{L^2}^2
+\|\mean(\q\cdot\w_n)\|_{L^2}^2.
\end{align*}
We estimate
\begin{align*}
\frac{1}{\cos(\theta+\tau)}
\Re\left(e^{-i(\theta+\tau)\sign\omega} \right.&\left.\spl \tilde K_n \u_n,\u_n\spr_{\Xn}\right)\\
&\geq \CtKn |\u_n|_{Y_n}^2 - \hmax \CYi \|\u_n\|_{\Xn}^2
- \CYii \|\u_n\|_{\Xn}|\u_n|_{Y_n}
\end{align*}
with constants $\CYi,\CYii>0$.
Thus
\begin{align*}
\frac{1}{\cos(\theta+\tau)}
\Re\Big(e^{-i(\theta+\tau)\sign\omega} & \spl B_n \u_n,\u_n\spr_{\Xn}\Big)\\
&\hspace{-5mm}\geq \frac{\CtB}{2}\|\u_n\|_{\Xn}^2+\CtKn |\u_n|_{Y_n}^2
-\hmax \CYi \|\u_n\|_{\Xn}^2 - \CYii \|\u_n\|_{\Xn}|\u_n|_{Y_n}\\
&\hspace{-5mm}\geq (\CtB/4-\hmax \CYi) \|\u_n\|_{\Xn}^2+(\CtKn-\CYii^2/\CtB^2) |\u_n|_{Y_n}^2.
\end{align*}
Now we choose $\CtKn>\CYii^2/\CtB^2$ and obtain the uniform stability of $B_n$ for large enough index $n$.
\\

\emph{4.\ step: asymptotic consistency of $B_n$.}\quad
Similar to the discrete setting it holds that $AT=B+K$ with
\begin{align*}
% v_n, v_n'
\nonumber
\spl B\u,&\u' \spr_{\IX}:=\\
&
\inner{c_s^2\rho\div \v,\div \v'}
-\inner{\rho i\conv\v,i\conv \v'}
+\inner{c_s^2\rho \PLz(\q\cdot \w),\PLz(\q\cdot \w')}\\
% v, w', w, v'
&
-\inner{\rho i\conv\v,\opd \w'}
+\inner{\rho \opd \w,i\conv\v'}\\
% w, w'
&
+\inner{\rho\opd \w,\opd \w'}
+\inner{\rho(i\gamma+\m) \w,\w'}\\
&
+\inner{\v,\v'}
+\CKeta \spl \Keta\v,\Keta\v' \spr_{\bV}
+\inner{\Kuni\w,\Kuni \w'}\\
\nonumber \\
&
+\CtKn (\spl\v,\v'\spr
+\spl \Keta\v,\Keta\v' \spr_{\bV}
+\spl\Kuni\w,\Kuni \w'\spr
+\spl \mean(\q\cdot\w), \mean(\q\cdot\w') \spr)\\
% v, v'
\nonumber\\
&
+\inner{c_s^2\rho\q\cdot\v,\div\v'}
+\inner{c_s^2\rho\div\v,\q\cdot\v'}
-\inner{\rho (\omega+i\angvel\times) \v, (\omega+i\angvel\times) \v'}\\
&
-\inner{\rho (\omega+i\angvel\times) \v, i\conv \v'}
-\inner{\rho i\conv \v, (\omega+i\angvel\times) \v'}
-i\omega\inner{\gamma\rho\v,\v'}
-\inner{\rho\m \v,\v'}\\
% v, w'
\nonumber\\
&
-\inner{\rho\m \v,\w'}
-i\omega\inner{\gamma\rho\v,\w'}
-\inner{c_s^2\rho \PLz(\q\cdot \v),\PLz(\q\cdot \w')}\\
&
-\inner{\rho (\omega+i\angvel\times) \v, \opd \w'}
-\inner{c_s^2\rho (\div+\PLz\q\cdot)\v,\Kuni\w'}
+\spl c_s^2\rho \mean(\q\cdot\v) , \div\w' \spr\\
% w, v'
\nonumber\\
&
+\inner{\rho\m \w,\v'}
+i\omega\inner{\gamma\rho\w,\v'}
+\inner{c_s^2\rho \PLz(\q\cdot \w),\PLz(\q\cdot \v')}\\
&
+\inner{\rho \opd \w, (\omega+i\angvel\times) \v'}
+\inner{c_s^2\rho \Kuni\w, (\div+\PLz\q\cdot)\v'}
-\spl c_s^2\rho \div\w, \mean(\q\cdot\v') \spr\\
% w, w'
\nonumber\\
&
-\spl c_s^2\rho \mean(\q\cdot\w) , \div\w' \spr
-\spl c_s^2\rho \div\w, \mean(\q\cdot\w') \spr
-\spl c_s^2\rho \Kuni\w,\Kuni\w' \spr 
\end{align*}
and
\begin{align*}
&\spl K\u,\u'\spr_{\IX}:=\\
&-\CtKn \spl\v,\v'\spr
-(\CKeta+\CtKn) \spl \Keta\v,\Keta\v' \spr_{\bV}
-(1+\CtKn) \spl\Kuni\w,\Kuni \w'\spr
-\CtKn \spl \mean(\q\cdot\w), \mean(\q\cdot\w') \spr
\end{align*}
for all $\u,\u'\in\IX$.
In addition the coercivity of $B$ follows along the same lines of the respective proof for $B_n$.
To prove the asymptotic consistency of $B, B_n$ we first show the asymptotic consistency of $K, K_n$.
Thus let $\u\in\IX$ be given.
We need to show that $\lim_{n\to\infty}\|p_nK\u-K_np_n\u\|_{\Xn}=0$.
Let $\u_n'\in\Xn$, $\|\u_n'\|_{\Xn}=1$, $n\in\IN$ be such that
$\|p_nK\u-K_np_n\u\|_{\Xn} \leq |\spl p_nK\u-K_np_n\u,\u_n'\spr_{\Xn}|+1/n$.
Let $\IN'\subset\IN$ be an arbitrary subsequence.
Due to \cref{lem:weakcompactness} there exist $\IN''\subset\IN'$ and $\u'\in\IX$ such that
$\u_n'\xrightharpoonup{\bL^2}\u'$, $ \div\u_n'\xrightharpoonup{L^2}\div\u'$, $ \Db\u_n'\xrightharpoonup{\bL^2}\conv\u'$, and we conveniently compute
\begin{align*}
\spl p_n K\u,\u_n'\spr_{\Xn}
&=\spl \div K\u,\div\u_n'\spr+\spl K\u,\u_n'\spr+\spl\conv K\u,\Db\u_n'\spr\\
&\xrightarrow{n\in\IN''} \spl \div K\u,\div\u'\spr+\spl K\u,\u'\spr+\spl\conv B\u,\conv\u'\spr
=\spl K\u,\u'\spr_{\IX}.
\end{align*}
On the other hand we estimate
\begin{align*}
|\spl \v,\v_n'\spr - \spl \PVn p_n\u,\v_n'\spr|
&=|\spl P_V\u - \pi_n^d \tPV p_n\u,\v_n'\spr|\\
&\lesssim |\spl P_V\u - \pi_n^d \tPV \u,\v_n'\spr| + \dn(\u,p_n\u)\\
&\lesssim |\spl P_V\u - \pi_n^d P_V \u,\v_n'\spr| + \dn(\u,p_n\u) + \|(\PLz-\pi_n^l)(\q\cdot\u)\|_{L^2}\\
&\lesssim \hmax \|P_V\u\|_{\bH^1} + \dn(\u,p_n\u) + \|(\PLz-\pi_n^l)(\q\cdot\u)\|_{L^2}\\
\end{align*}
and
\begin{align*}
|\spl \Keta\v,\Keta\tPV\u_n' \spr_{\bV} &- \spl \Keta\tPV p_n\u,\Keta\tPV\u_n' \spr_{\bV}|\\
&=|\spl \Keta (P_V\u-\tPV p_n \u),\Keta\tPV\u_n' \spr_{\bV}|\\
&\lesssim |\spl \Keta (P_V\u-\tPV \u),\Keta\tPV\u_n' \spr_{\bV}| + \dn(\u,p_n\u)\\
&\lesssim \|(\PLz-\pi_n^l)(\q\cdot\u)\|_{L^2} + \dn(\u,p_n\u)
\end{align*}
and
\begin{align*}
|\spl \Kuni \w,\Kuni \w_n' \spr_{\bV}&-\spl \Kuni \w_n(p_n\u),\Kuni \w_n' \spr|\\
&=|\spl \Kuni (\w-\w_n(p_n\u)),\Kuni \w_n' \spr|\\
&=|\spl \Kuni (\u-P_V\u-(p_n\u-\PVn p_n\u)),\Kuni \w_n' \spr|\\
&\lesssim\|\u-P_V\u-(p_n\u-\PVn p_n\u)\|_{H(\div)}\\
&\lesssim\|P_V\u-\PVn p_n\u\|_{H(\div)}+\dn(\u,p_n\u)\\
&\lesssim\|p_n P_V\u-\PVn p_n\u\|_{H(\div)}+\dn(\u,p_n\u)+\dn(P_V\u,p_n P_V\u)
\end{align*}
and
\begin{align*}
|\spl \mean(\q\cdot\w),\mean(\q\cdot\w_n')\spr_{\bV}&-\spl \mean(\q\cdot\w_n(p_n\u)),\mean(\q\cdot\w_n')\spr|\\
&\lesssim\|\u-P_V\u-(p_n\u-\PVn p_n\u)\|_{H(\div)}\\
&\lesssim\|p_n P_V\u-\PVn p_n\u\|_{H(\div)}+\dn(\u,p_n\u)+\dn(P_V\u,p_n P_V\u).
\end{align*}
Thus
\begin{align*}
\lim_{n\to\infty}
|\spl K_n p_n\u,\u_n'\spr_{\Xn}&+\CtKn \spl\v,\v_n'\spr
+(\CKeta+\CtKn) \spl \Keta\v,\Keta\tPV\u_n' \spr_{\bV}
+(1+\CtKn) \spl\Kuni\w,\Kuni\w_n'\spr\\
&+\CtKn \spl \mean(\q\cdot\w), \mean(\q\cdot\w_n') \spr|
=0.
\end{align*}
We further use the operator $S:=\nabla\big((\div+\PLz\q\cdot+\Kuni)\nabla\big)^{-1}\in L(L^2_0,\bV)$ (whereat we changed the space $\bL^2$ to $\bV$ compared to the proof of \cref{lem:Kncompact}) and compute
\begin{align*}
\spl \Keta\v,\Keta\tPV\u_n' \spr_{\bV}
&=\spl \Keta^*\Keta\v,\tPV\u_n' \spr_{\bV}\\
&=\spl \Keta^*\Keta\v,P_V\u_n' \spr_{\bV}
+\spl \Keta^*\Keta\v, S (\pi_n^l-\PLz)(\q\cdot\u_n') \spr_{\bV}\\
&=\spl P_V^*\Keta^*\Keta\v,\u_n' \spr_{H_0(\div)}
+\spl (\pi_n^l-\PLz)S^*\Keta^*\Keta\v, \q\cdot\u_n' \spr_{L^2}\\
&\xrightarrow{n\in\IN''} \spl P_V^*\Keta^*\Keta\v,\u' \spr_{H_0(\div)}
=\spl \Keta\v,\Keta\v' \spr_{\bV}
\end{align*}
and hence
\begin{align*}
&\lim_{n\in\IN''} \Big(
\CtKn \spl\v,\v_n'\spr
+(\CKeta+\CtKn) \spl \Keta\v,\Keta\tPV\u_n' \spr_{\bV}
+(1+\CtKn) \spl\Kuni\w,\Kuni\w_n'\spr\\
&\qquad+\CtKn \spl \mean(\q\cdot\w), \mean(\q\cdot\w_n') \spr \big)\\
&=\CtKn \spl\v,\v'\spr
+(\CKeta+\CtKn) \spl \Keta\v,\Keta\v' \spr_{\bV}
+(1+\CtKn) \spl\Kuni\w,\Kuni\w'\spr
+\CtKn \spl \mean(\q\cdot\w), \mean(\q\cdot\w') \spr\\
&=\spl K\u,\u'\spr_{\IX}.
\end{align*}
Thus we obtained that $\lim_{n\to\infty} \|p_n K\u-K_n p_n \u\|_{\Xn}$ for each $\u\in\IX$.
It remains to recall $B_n=A_nT_n-K_n$, $B=AT-K$ and to estimate
\begin{align*}
\|(p_nB-B_n&p_n)\u\|_{\Xn}\\
&\leq \|(p_nK-K_np_n)\u\|_{\Xn} + \|(p_nAT-A_nT_np_n)\u\|_{\Xn}\\
&\leq \|(p_nK-K_np_n)\u\|_{\Xn} + \|(p_nA-A_np_n)T\u\|_{\Xn}
+ \|A_n\|_{L(\Xn)} \|(p_nT-T_np_n)\u\|_{\Xn}.
\end{align*}
Thus $\lim_{n\to\infty}\|(p_nB-B_np_n)\u\|_{\Xn}=0$ follows from the just proven asymptotic consistency of $K_n, K$, \cref{thm:asympcons}, \cref{lem:Tnasympcons} and from the uniform boundedness of $(A_n)_{n\in\IN}$.
\end{proof}

\subsection{Convergence results}

\begin{theorem}\label{thm:convergence}
Let $\bf\in\bL^2$ and $\u\in\IX$ be the solution to $a(\u,\u')=\spl\bf,\u'\spr$ for all $\u'\in\IX$.
Then there exists $n_0>0$ such that for all $n>n_0$ the solution $\u_n\in\Xn$ to $\sesn(\u_n,\u_n')=\spl\bf,\u_n'\spr$ for all $\u_n'\in\Xn$ exists and $\lim_{n\to\infty} \dn(\u,\u_n)=0$.
In addition, if $\u\in\bH^{2+s}$, $\rho\in W^{1+s,\infty}$ and $\vel\in\bW^{1+s,\infty}$ with $s>0$, then $\dn(\u,\u_n) \lesssim \hmax^{\min(1+s,k)}+\hmax^{\min(s,l)}$.
\end{theorem}
\begin{proof}
Due to \cref{thm:asympcons,thm:dis-wTc} we can apply \cref{thm:Tcomp}.
Since $A$ is injective and hence bijective \cref{lem:DASstable} yields that $(A_n)_{n\in\IN}$ is stable.
Let $\g\in\IX$ be such that $\spl\g,\u'\spr_{\IX}=\spl\bf,\u'\spr_{\bL^2}$ for all $\u'\in\IX$ and $\g_n\in\Xn$ be such that $\spl\g_n,\u'_n\spr_{\Xn}=\spl\bf,\u'_n\spr_{\bL^2}$ for all $\u'_n\in\Xn$.
To obtain that $\lim_{n\to\infty} \dn(\u,\u_n)=0$ it remains to show that $\lim_{n\to\infty}\|p_n\g-\g_n\|_{\Xn}=0$.
We proceed conveniently and choose $\u_n'\in\Xn, \|\u_n'\|_{\Xn}=1, n\in\IN$ such that $\|p_n\bg-\bg_n\|_{\Xn}\leq |\spl p_n\g-\g_n,\u_n'\spr_{\Xn}|+1/n$.
For an arbitrary subsequence $\IN'\subset\IN$ we choose $\u'\in\IX$ and $\IN''\subset\IN'$ as in \cref{lem:weakcompactness} and obtain that
\begin{align*}
\spl p_n\bg-\bg_n,\u_n'\spr_{\Xn}
=\spl \div\g,\div\u_n'\spr+\spl \g,\u_n'\spr+\spl \conv\g,\Db\u_n'\spr -\spl\bf,\u_n'\spr
\xrightarrow{n\in\IN''} \spl \g,\u'\spr_{\IX}-\spl\bf,\u'\spr=0,
\end{align*}
from which it follows that $\lim_{n\to\infty}\|p_n\g-\g_n\|_{\Xn}=0$, and hence $\lim_{n\to\infty} \dn(\u,\u_n)=0$.
To obtain the convergence rate we first estimate
\begin{align*}
\dn(\u,\u_n) \leq \dn(\u,p_n\u)+\|p_n\u-\u_n\|_{\Xn}
\lesssim \dn(\u,p_n\u)+\|A_n(p_n\u-\u_n)\|_{\Xn}
\end{align*}
and further compute that
\begin{align*}
\|A_n(p_n\u-\u_n)\|_{\Xn}
&=\sup_{\|\u_n'\|_{\Xn}=1} |\sesn(p_n\u-\u_n,\u_n')|
=O\big(\dn(\u,p_n\u), n\to\infty\big)\\
&+\sup_{\|\u_n'\|_{\Xn}=1} \big|\spl c_s^2\rho\div \u,\div \u_n'\spr
-\spl\rho\opd \u,(\omega+i\Db+i\angvel\times) \u_n'\spr\\
&+\spl\div\u,\nabla p\cdot \u_n'\spr
+\spl\nabla p\cdot\u,\div \u_n'\spr
+\spl(\hess(p)-\rho\hess(\phi)) \u,\u_n'\spr\\
&-i\omega\spl\gamma\rho\u,\u_n'\spr-\spl \bf,\u_n'\spr_{\bL^2}\big|.
\end{align*}
The next step is to integrate by parts
\begin{align*}
\spl c_s^2\rho\div \u,\div \u_n'\spr=-\spl \nabla(c_s^2\rho\div \u),\u_n'\spr,\qquad
\spl\nabla p\cdot\u,\div \u_n'\spr=-\spl\nabla(\nabla p\cdot\u),\u_n'\spr,
\end{align*}
whereat the left arguments in the $\bL^2$-scalar products of the former right hand-sides are in $\bL^2$ due to the assumptions $\u\in\bH^2$, $c_s, \rho\in W^{1,\infty}$ and $p\in W^{2,\infty}$.
The integration by parts of the operator $\Db$ is performed similarly as in the proof of \cref{lem:weakcompactness}.
Thus let $\bpsi_n\in\bQ_n$ be a suitable $\bH^1$ projection of $\rho\opd\bu$, e.g.\ $\bpsi_n=\calJ_n(\rho\opd\u)$ with $\calJ_n$ as in \cite[\bb6.4)]{ErnGuermond:16}
We compute
\begin{align*}
\spl\rho(\omega+i\conv&+i\angvel\times)\u,\Db\u_n'\spr
=\spl\bpsi_n,\Db\u_n'\spr
+\spl\rho\opd\u-\bpsi_n,\Db\u_n'\spr
\end{align*}
and
\begin{align*}
\spl\bpsi_n,\Db\u_n'\spr
&=\sum_{\tet\in\calT_n} \spl\bpsi_n,\conv\u_n'+\LO\u_n'\spr_{\bL^2(\tet)}
=\sum_{\tet\in\calT_n} \spl\bpsi_n,\conv\u_n'\spr_{\bL^2(\tet)}
-\spl\avg{\bpsi_n},\bjump{\u_n'}\spr_\Find\\
&=\sum_{\tet\in\calT_n} \spl\bpsi_n,\conv\u_n'\spr_{\bL^2(\tet)}
-\spl\bpsi_n,(\nv\cdot\vel)\u_n'\spr_{\bL^2(\partial\tet)}\\
&=-\spl(\conv+\div(\vel))\bpsi_n,\u_n'\spr\\
&=-\spl(\conv+\div(\vel))\rho(\omega+i\conv+i\angvel\times)\u,\u_n'\spr\\
&+\spl(\conv+\div(\vel))\big(\bpsi_n-\rho(\omega+i\conv+i\angvel\times)\u\big),\u_n'\spr\\
&=-\spl\rho\conv(\omega+i\conv+i\angvel\times)\u,\u_n'\spr
-\spl\div(\rho\vel)(\omega+i\conv+i\angvel\times)\u,\u_n'\spr\\
&+\spl(\conv+\div(\vel))\big(\bpsi_n-\rho(\omega+i\conv+i\angvel\times)\u\big),\u_n'\spr.
\end{align*}
Thus
\begin{align*}
\sup_{\|\u_n'\|_{\Xn}=1} \big|&\spl c_s^2\rho\div \u,\div \u_n'\spr
-\spl\rho\opd \u,(\omega+i\Db+i\angvel\times) \u_n'\spr
+\spl\div\u,\nabla p\cdot \u_n'\spr\\
&+\spl\nabla p\cdot\u,\div \u_n'\spr
+\spl(\hess(p)-\rho\hess(\phi)) \u,\u_n'\spr
-i\omega\spl\gamma\rho\u,\u_n'\spr-\spl \bf,\u_n'\spr_{\bL^2}\big|\\
&\lesssim \|\rho(\omega+i\conv+i\angvel\times)\u-\bpsi_n\|_{\bH^1}
\lesssim \hmax^{\min(s,l)}
\end{align*}
due to the properties of $\calJ_n$ and \eqref{eq:EstApr}.
Since we can estimate $\dn(\u,p_n\u)$ with \cref{lem:dnpihn} the claim follows.
\end{proof}
\begin{remark}
The unusual regularity assumptions and convergence rate $\hmax^{\min(1+s,k)}+\hmax^{\min(s,l)}$ of \cref{thm:convergence} deserve some discussion.
First we note that for a right hand-side $\bf\in\bL^2$ it follows that $-\nabla(c_s^2\rho\div\u)+\rho\conv\conv\u \in\bL^2$, which although does not allow us to deduce neither $\nabla(c_s^2\rho\div\u)\in\bL^2$ nor $\rho\conv\conv\u \in\bL^2$.
This is the reason why we explicitly need to assume that $\u\in\bH^2$.
If the polynomial degree $l$ of the lifting operator $\LO$ is choosen as $l=k$, and if $k\leq s$, then we obtain the convenient rate $\hmax^k$.
However, if the solution and the parameters have only a maximal regularity $s<\infty$ and $k\geq 1+s$, then we obtain only the rate $\hmax^s$ which is one power less then the conveniently expected rate.
The reason for this unusual result is that we employ a DG method without a ($\vel$-jump) stabilization term, which results in a weaker norm.
\end{remark}
\begin{remark}
In principle the previous analysis of this article can also be applied to other DG variants such as the symmetric interior penalty method, whereat we note that the sign of the coefficient of the penalty term $\spl \h^{-1} \bjump{\u_n},\bjump{\u_n'} \spr_\Find$ in the sesquilinear form needs to be negative (or have a suitable complex sign).
Of course, the norm $\|\cdot\|_{\Xn}$ and everything related would need to be adapted.
Further note that it follows along the lines of the proof of \cite[Theorem~8]{AHLS:22} that there exists a constant $C>0$ such that
\begin{align*}
\|(\h\rho)^{\sfrac12}\avg{\conv\u_n}\|_\Find
&\leq C (\|\rho^{\sfrac12}\conv\u_n\|_{\bL^2}+\|\rho^{\sfrac12}\u_n\|_{\bL^2}),\\
\|(\h\rho)^{\sfrac12}\avg{\opd\u_n}\|_\Find
&\leq C (\|\rho^{\sfrac12}\opd\u_n\|_{\bL^2}+\|\rho^{\sfrac12}\u_n\|_{\bL^2})
\end{align*}
for all $\u_n\in\Xn$, $n\in\IN$, whereat $\conv\bu_n$ is interpreted piece-wise with respect to the mesh $\calT_n$.
Remark the unconvenient term $\|\rho^{\sfrac12}\u_n\|_{\bL^2}$ in the former right hand-sides.
Thence it follows that the for a large enough penalty parameter $\ALpen>0$ there holds the coercivity estimate
\begin{align}\label{eq:EstSIP}
\begin{aligned}
\big|\|\rho^{\sfrac12}(i\omega+&i\conv+i\angvel\times)\w_n\|_{\bL^2}^2
+\ALpen\|\h^{-\sfrac12}\bjump{\w_n}\|_\Find^2
+i\omega\|(\rho\gamma)^{\sfrac12}\w_n\|_{\bL^2}^2\\
&-\spl \rho\avg{\opd\w_n},\bjump{\w_n} \spr_\Find
-\spl \bjump{\w_n},\rho\avg{\opd\w_n} \spr_\Find\big|\\
&\gtrsim \|\opd\w_n\|_{\bL^2}^2
+\|\h^{-\sfrac12}\bjump{\w_n}\|_\Find^2
+\|(\rho\gamma)^{\sfrac12}\w_n\|_{\bL^2}^2.
\end{aligned}
\end{align}
Now a crux is that in proof of \cref{thm:dis-wTc} we multiply with $e^{-i(\theta+\tau)\sign\omega}$ and hence to real part of the coefficient of $\|(\rho\gamma)^{\sfrac12}\u_n\|_{\bL^2}^2$ depends on $\theta$ and will become small for large $\theta$.
Thus $\ALpen$ will depend on $\theta$ as well and has to be chosen sufficiently large to guarantee a coercivity estimate of kind \eqref{eq:EstSIP}.
Until now this produces no severe drawbacks.
However, when repeating the respective estimates in the proof of \cref{thm:dis-wTc} the additional term $-\ALpen\|\h^{-\sfrac12}\bjump{\pi_n^d\tilde\v}\|_\Find^2$ needs to be estimated by $|\nabla\tilde\v|_{\bH^1_{c_s^\rho}}^2$, which leads to a more restrictive assumption on the smallness of the Mach number.
We conclude that the penalty parameter $\ALpen$ needs to be balanced in a nontrivial way to guarantee a coercivity estimate for $\w_n$ while avoiding an unnecessary confining assumption on the smallness of the Mach number.
\end{remark}

% \section{Computational examples}\label{sec:numex}
% 
% todo
% 
% \subsection{Implementation}
% 
% The assembling of the system matrix associated to $\sesn(\cdot,\cdot)$ requires the realization of the operator $\LO$.
% Note that $\LO=\sum_{\face\in\calF_n^\mathrm{int}} \LO^\face$ and each $\LO^\face$ is a \emph{local nonlocal operator}.
% This phrase is a bit confusing, so let us elaborate.
% Since $\LO^\face$ is defined as a solution to a variational problem it is a nonlocal operator.
% However, all functions which have a support contained in $\ol{\Omega\setminus{\ol{\tet_1\cup\tet_2}}}$ are in the kernel of $\LO^\face$ and in addition it holds that $\supp\LO^\face\u_n\subset \ol{\tet_1\cup\tet_2}$ for each $\u_n\in\Xn$, which is a localness property.
% We remark that the lifting operator leads to an increased sparsity pattern.
% Indeed, all basis functions which are associated to a $\tet\in\calT_n$ or $\face\in\calF(\tet)$ couple with the respective functions associated to all neighboring elements $\tilde\tet$ of all neighbors of $\tet$.
% To get a more local operator we introduce a hybrid formulation.
% Thus we introduce a facet space $Y_n\subset \bL^2(\calF_n)$ of order $k$.
% 
% \subsection{Robustness w.r.t.\ $\rho,c_s$}

\bibliographystyle{siam}
\bibliography{bib}

\end{document}